# PERTURBATION SELECTION AND INFLUENCE MEASURES IN LOCAL INFLUENCE ANALYSIS

By Hongtu Zhu,[1] Joseph G. Ibrahim,[2] Sikyum Lee
and Heping Zhang[3]

*University of North Carolina at Chapel Hill, University of North Carolina at Chapel Hill, The Chinese University of Hong Kong and Yale University*

Cook's [*J. Roy. Statist. Soc. Ser. B* **48** (1986) 133–169] local influence approach based on normal curvature is an important diagnostic tool for assessing local influence of minor perturbations to a statistical model. However, no rigorous approach has been developed to address two fundamental issues: the selection of an appropriate perturbation and the development of influence measures for objective functions at a point with a nonzero first derivative. The aim of this paper is to develop a differential–geometrical framework of a perturbation model (called the perturbation manifold) and utilize associated metric tensor and affine curvatures to resolve these issues. We will show that the metric tensor of the perturbation manifold provides important information about selecting an appropriate perturbation of a model. Moreover, we will introduce new influence measures that are applicable to objective functions at any point. Examples including linear regression models and linear mixed models are examined to demonstrate the effectiveness of using new influence measures for the identification of influential observations.

**1. Introduction.** Assessing local influence of perturbing a statistical model has been an active area of statistical research in the past twenty years since the seminal work of Cook [7]. See, for example, Beckman, Nachtsheim and Cook [3], Tsai and Wu [26], St. Laurent and Cook [25], Wu and Luo [31, 32], Ouwens, Tan and Berger [21], Pan and Fang [22] and Zhu and Lee [38],

Received May 2006; revised September 2006.
[1]Supported by NSF Grant SES-06-43663.
[2]Supported by NIH Grants GM 70335 and CA 74015.
[3]Supported by NIDA Grants DA016750 and DA017713.

among many others. The key idea of the local influence approach is to utilize the concept of normal curvature in differential geometry (Efron [10] and Bates and Watts [2]) in assessing the local behavior of the likelihood displacement function. Zhu and Lee [37] proposed a generalization of Cook's [7] approach based on the Q-function in the EM algorithm to assess local influence of a small perturbation to a class of models with incomplete data. Zhu, He and Fung [41] developed a local influence method for generalized partial linear models for longitudinal data. Lee and Tang [17] examined local influence in structural equation models. Zhu and Zhang [40] used measures of local influence to assess the extent of discrepancy between a hypothetical model and the underlying model from which the data are generated. The local influence approach is also useful for sensitivity analysis of missing data modeling (Verbeke and Molenberghs [27] and Verbeke et al. [28]).

The aim of this paper is to construct influence measures in assessing local influence of perturbations to a statistical model. Specifically, we address two important issues related to the local influence approach: the appropriate choice of a perturbation vector and development of influence measures for assessing an objective function at any point.

The first issue, selecting an appropriate perturbation vector, has been largely neglected. This issue, however, is central to the development of the local influence approach, because arbitrarily perturbing a model may lead to inappropriate inference about the cause (e.g., influential observations) of a large effect. For instance, when data form clusters (e.g., subjects in longitudinal studies and families in genetic studies), then perturbing a cluster with more observations likely produces a larger effect. However, to the best of our knowledge, no diagnostic methods have ever been developed to take into account the differing number of observations in each cluster (see Section 4 for further discussion). Moreover, because the components of a perturbation vector may not be orthogonal to each other, special care should be taken when we interpret local influence measures from such a perturbation. Thus, it is desirable to measure the amount of perturbation, the extent to which each component of a perturbation vector contributes to, and the degree of orthogonality for the components of a perturbation vector.

The second issue is the development of influence measures for objective functions at a point with a nonzero first derivative. Fung and Kwan [11] showed that the normal curvature is not scale invariant and provided some examples to illustrate that ambiguous conclusions may be drawn when applied to objective functions with a nonzero first derivative. However, they did not provide any methodology to address this drawback. Conformal normal curvature (Poon and Poon [23]) is invariant under the conformal reparametrization at a point with a zero first derivative (called a critical point), but it is not scale invariant at a point with a nonzero first derivative. These difficulties have limited the application of normal curvature to objective functions that have zero first derivative at the critical point.

We will introduce a geometrical structure, called the perturbation manifold, for a perturbation model and use its associated metric tensor and affine connection to address the above two issues. The metric tensor of the perturbation manifold can measure the amount of perturbation and the orthogonality between the different components of a perturbation vector. Thus, the properties of metric tensors (e.g., positive definiteness) can be used to choose an appropriate perturbation to a statistical model. Furthermore, once an appropriate perturbation is chosen, we use the first and second derivatives of the objective function (e.g., the likelihood displacement function in Cook [7]) to construct influence measures. These influence measures can be easily applied to any objective function evaluated at any point in order to quantify the local influence of minor perturbations to a statistical model.

## 2. Perturbation manifold and influence measures.

2.1. *Motivation.* Let $p(\mathbf{Y}|\boldsymbol{\theta})$ be the probability function for an $M(n) \times 1$ random vector $\mathbf{Y}^T = (Y_1^T, \ldots, Y_n^T)$, parameterized by an unknown parameter vector $\boldsymbol{\theta} = (\theta_1, \ldots, \theta_q)^T$ in an open subset $\Theta$ of $R^q$. In addition, each $Y_i$ is an $m_i \times 1$ random vector, where $\sum_{i=1}^n m_i = M(n)$. For instance, in longitudinal studies, $m_i$ may represent the number of observations in the $i$th cluster. On the basis of the assumed model $p(\mathbf{Y}|\boldsymbol{\theta})$ and observations in $\mathbf{y}^T = (\mathbf{y}_1^T, \ldots, \mathbf{y}_n^T)$, we can then carry out statistical inference, such as estimation and hypothesis testing.

Let $\boldsymbol{\omega} = (\omega_1, \ldots, \omega_p)^T$ be a perturbation vector and let $\boldsymbol{\omega}$ vary in $\Omega \subset R^p$. If a perturbation vector $\boldsymbol{\omega}$, which is introduced to perturb $p(\mathbf{Y}|\boldsymbol{\theta})$, has a large effect, then it is important to know the cause (e.g., influential observations or invalid model assumptions) of such a large effect. Therefore, it is important to develop statistical methods to quantify the effect of perturbing a statistical model and pinpoint the potential cause.

In Cook [7] a general method was developed to assess the local influence of perturbing a statistical model by introducing $\boldsymbol{\omega}$ into $p(\mathbf{Y}|\boldsymbol{\theta})$, denoted by $p(\mathbf{Y}|\boldsymbol{\theta}, \boldsymbol{\omega})$. The proposed methodology is based on the directional curvature of an influence graph, which is defined as

$$(1) \qquad \mathrm{IG}(\boldsymbol{\omega}) = (\boldsymbol{\omega}^T, f(\boldsymbol{\omega}))^T,$$

where $f : R^p \to R^1$ is a sufficiently smooth (differentiable a certain number of times) objective function. Consider the straight line $\boldsymbol{\omega}(t)$: $\boldsymbol{\omega}(t) = \boldsymbol{\omega}^0 + t\mathbf{h}$ in Euclidean space $R^p$ and the lifted line $\mathrm{IG}_{\mathbf{h}}(\boldsymbol{\omega}(t))$ for any nonzero vector $\mathbf{h}$, where $\boldsymbol{\omega}^0$ is a fixed column vector in $R^p$. The tangent vector and upward-point normal vector of the lifted line are, respectively, given by $\binom{\mathbf{I}_{p;}}{\nabla_f^T}$ and $(1 + \nabla_f^T \nabla_f)^{-1/2} \binom{-\nabla_f}{1}$, where $\nabla_f = (\partial f(\boldsymbol{\omega})/\partial \omega_i)$ is evaluated at $\boldsymbol{\omega}^0$ and $\mathbf{I}_p$ is the $p \times p$ identity matrix. The normal curvature of the influence graph (Cook [7]) is given by

$$(2) \qquad C_{\mathbf{h}} = \frac{1}{(1 + \nabla_f^T \nabla_f)^{1/2}} \frac{\mathbf{h}^T H_f \mathbf{h}}{\mathbf{h}^T (\mathbf{I}_p + \nabla_f \nabla_f^T) \mathbf{h}},$$

where $H_f$ denotes the matrix $(\partial^2 f(\omega)/\partial \omega_i \, \partial \omega_j)$ evaluated at $\omega^0$. The maximum value of $C_{\mathbf{h}}$ and the corresponding direction have been widely used to assess the effects of using $\omega(t) = \omega^0 + t\mathbf{h}$ to perturb a statistical model. Poon and Poon [23] defined a conformal normal curvature at $\omega^0$ in the direction $\mathbf{h}$ as

$$(3) \qquad B_{\mathbf{h}} = \frac{1}{\|H_f\|_M} \frac{\mathbf{h}^T H_f \mathbf{h}}{\mathbf{h}^T (\mathbf{I}_p + \nabla_f \nabla_f^T) \mathbf{h}},$$

where $\| \cdot \|_M$ denotes the norm of a matrix such that $\|H_f\|_M = \sqrt{\mathrm{tr}[H_f]^2}$. However, $C_{\mathbf{h}}$ is not scale invariant at any $\omega$ with $\nabla_f \neq 0$, because the normal curvature of $\widehat{\mathrm{IG}}(\omega) = (\omega^T, k f(\omega))^T$ is given by

$$\widehat{C_{\mathbf{h}}} = \frac{1}{(1 + k^2 \nabla_f^T \nabla_f)^{1/2}} \frac{k \mathbf{h}^T H_f \mathbf{h}}{\mathbf{h}^T (\mathbf{I}_p + k^2 \nabla_f \nabla_f^T) \mathbf{h}} \neq C_{\mathbf{h}},$$

which may lead to ambiguous conclusions as $k$ varies (Fung and Kwan [11]). The same problem also arises with $B_{\mathbf{h}}$, that is, $B_{\mathbf{h}}$ is not scale invariant. In particular, Fung and Kwan [11] argued that the conclusions drawn from the new graph $\widehat{\mathrm{IG}}(\omega)$ should be the same as those from the old graph $\mathrm{IG}(\omega) = (\omega^T, f(\omega))^T$.

2.2. *Statistical perturbation manifold.* We use $p(\mathbf{Y}|\boldsymbol{\theta}, \omega)$ to denote the density function such that $\int p(\mathbf{Y}|\boldsymbol{\theta}, \omega) \, d\mathbf{Y} = 1$. To assess the local influence of a model perturbation, we are primarily interested in the behavior of $p(\mathbf{Y}|\boldsymbol{\theta}, \omega)$ as a function of $\omega$ around $\omega^0$, not the parameter vector $\boldsymbol{\theta}$. From now on, $\boldsymbol{\theta}$ is assumed to be known or be fixed at a given value (e.g., the maximum likelihood estimate) and $p(\mathbf{Y}|\boldsymbol{\theta}, \omega^0) = p(\mathbf{Y}|\boldsymbol{\theta})$. Moreover, $p(\mathbf{Y}|\boldsymbol{\theta}, \omega)$ satisfies the four regularity conditions on page 16 of Amari [1] and $\omega^0$ represents no perturbation.

The perturbed model $p(\mathbf{Y}|\boldsymbol{\theta}, \omega)$ is characterized by a set of perturbations $\omega$, which has a natural geometrical structure (Amari [1]). The perturbed model $M = \{p(\mathbf{Y}|\boldsymbol{\theta}, \omega) : \omega \in \Omega\}$ can be regarded as a $p$-dimensional manifold. When a coordinate system $\omega$ is given, $\mathbf{e}_i$ $(i = 1, \dots, p)$ are the natural basis of the tangent space $T_{\omega}$ of $M$ associated with the coordinate system (Amari [1] and Li and McCullagh [18]). Let $T_{\omega}^{(1)}$ be the vector space of $M$ at $\omega$, which is spanned by $p$ functions $\partial_i \ell(\omega|\mathbf{Y}, \boldsymbol{\theta})$, where $\ell(\omega|\mathbf{Y}, \boldsymbol{\theta}) = \log p(\mathbf{Y}|\boldsymbol{\theta}, \omega)$. A natural isomorphism exists between these two tangent vector spaces $T_{\omega}$ and

$T_{\boldsymbol{\omega}}^{(1)}$. The vector space $T_{\boldsymbol{\omega}}^{(1)}$ is called the 1-representation of the tangent space of $M$. For any tangent vector $\mathbf{h} = \sum_{i=1}^{p} h^i \mathbf{e}_i \in T_{\boldsymbol{\omega}}$, the 1-representation $\mathbf{h}(\mathbf{Y})$ of $\mathbf{h}$ in $T_{\boldsymbol{\omega}}^{(1)}$ is given by $\mathbf{h}(\mathbf{Y}) = \sum_{i=1}^{p} h^i \, \partial_i \ell(\boldsymbol{\omega}|\mathbf{Y}, \boldsymbol{\theta})$, where $\partial_i = \partial/\partial\omega_i$.

DEFINITION 1. The inner product of two basis operators $\partial_i$ and $\partial_j$ is

$$(4) \qquad g_{ij}(\boldsymbol{\omega}) = \langle \partial_i, \partial_j \rangle = \mathrm{E}_{\boldsymbol{\omega}}[\partial_i \ell(\boldsymbol{\omega}|\mathbf{Y}, \boldsymbol{\theta}) \, \partial_j \ell(\boldsymbol{\omega}|\mathbf{Y}, \boldsymbol{\theta})],$$

where $\mathrm{E}_{\boldsymbol{\omega}}$ denotes the expectation taken with respect to $p(\mathbf{Y}|\boldsymbol{\theta}, \boldsymbol{\omega})$. The $p^2$ quantities $g_{ij}(\boldsymbol{\omega})$, $i, j = 1, \ldots, p$, form the *metric* tensor.

The metric matrix $G(\boldsymbol{\omega}) = (g_{ij}(\boldsymbol{\omega}))$ is an expected Fisher information matrix with respect to the perturbation vector $\boldsymbol{\omega}$. The elements of $G(\boldsymbol{\omega})$ measure the amount of perturbation that all components of a perturbation vector $\boldsymbol{\omega}$ contribute to a statistical model. The $(i, i)$th element $g_{ii}(\boldsymbol{\omega})$ itself indicates the amount of perturbation introduced by the $i$th component of $\boldsymbol{\omega}$. The off-diagonal elements of $G(\boldsymbol{\omega})$ represent the association between different components of $\boldsymbol{\omega}$. For instance, let $r_{ij}(\boldsymbol{\omega}) = g_{ij}(\boldsymbol{\omega})/\sqrt{g_{ii}(\boldsymbol{\omega})g_{jj}(\boldsymbol{\omega})}$. A large absolute value of $r_{ij}(\boldsymbol{\omega})$ indicates strong association between the $i$th and $j$th components of $\boldsymbol{\omega}$. In particular, if $G(\boldsymbol{\omega})$ is a diagonal matrix, then all components of $\boldsymbol{\omega}$ are orthogonal to each other in the perturbed model (Cox and Reid [9]). Moreover, if $G(\boldsymbol{\omega})$ is not positive definite for a perturbation scheme, then $p$ operators $\partial_i$ are linearly dependent. Thus, this indicates that some components of the perturbation vector are redundant and these redundant components should be removed; for further discussion, see Section 3.3.2.

Based on the above discussion, an *appropriate perturbation* to a statistical model should satisfy at least two conditions as follows:

(a) $G(\boldsymbol{\omega})$ is positive definite in a small neighborhood of $\boldsymbol{\omega}^0$;
(b) the off-diagonal elements of $G(\boldsymbol{\omega})$ at $\boldsymbol{\omega}^0$ should be as small as possible.

Condition (a) is required to avoid any redundant components of $\boldsymbol{\omega}$. Condition (b) is required to ensure that we can easily pinpoint the cause of a large effect. For instance, if differing components of $\boldsymbol{\omega}$ are highly associated, then it is difficult to infer whether a large effect is caused by a single component of $\boldsymbol{\omega}$ or by several components of $\boldsymbol{\omega}$. Therefore, an appropriate perturbation requires that $G(\boldsymbol{\omega}^0)$ should be $\mathrm{diag}(g_{11}(\boldsymbol{\omega}^0), \ldots, g_{pp}(\boldsymbol{\omega}^0))$. Moreover, we can always choose a new perturbation vector $\tilde{\boldsymbol{\omega}}$, defined by

$$(5) \qquad \tilde{\boldsymbol{\omega}} = \boldsymbol{\omega}^0 + c^{-1/2} G(\boldsymbol{\omega}^0)^{1/2}(\boldsymbol{\omega} - \boldsymbol{\omega}^0),$$

such that $G(\tilde{\boldsymbol{\omega}})$ evaluated at $\boldsymbol{\omega}^0$ equals $c\mathbf{I}_p$, where $c > 0$. Therefore, without loss of generality, we assume that an appropriate perturbation $\boldsymbol{\omega}$ satisfies $G(\boldsymbol{\omega}^0) = c\mathbf{I}_p$. However, it is not generally possible to find a perturbation vector such that $G(\boldsymbol{\omega}) = c\mathbf{I}_p$ for all $\boldsymbol{\omega} \in \Omega$.

We introduce the following geometrical quantities for the perturbed model $M$ based on the metric tensor. First, the length $\|\mathbf{h}\|^2$ of a tangent vector $\mathbf{h} \in T_{\boldsymbol{\omega}}$ is given by

$$\|\mathbf{h}\|^2 = \langle \mathbf{h}, \mathbf{h} \rangle = \sum_{i,j} h^i h^j g_{ij}(\boldsymbol{\omega}) = \mathbf{h}^T G(\boldsymbol{\omega})\mathbf{h}. \tag{6}$$

Let $C : \boldsymbol{\omega}(t) = (\omega_1(t), \ldots, \omega_p(t))$ be a smooth curve on the manifold $M$ connecting two points $\boldsymbol{\omega}^1 = \boldsymbol{\omega}(t_1)$ and $\boldsymbol{\omega}^2 = \boldsymbol{\omega}(t_2)$. The distance $S(\boldsymbol{\omega}^1, \boldsymbol{\omega}^2)$ from $\boldsymbol{\omega}^1$ to $\boldsymbol{\omega}^2$ along the curve $C$ is given by

$$S(\boldsymbol{\omega}^1, \boldsymbol{\omega}^2) = \int_{t_1}^{t_2} \sqrt{\sum_{i,j} g_{ij}(\boldsymbol{\omega}(t)) \frac{d\omega_i(t)}{dt} \frac{d\omega_j(t)}{dt}} \, dt. \tag{7}$$

The skewness tensor $T$ and a family of affine connections $\Gamma^\alpha$ for any $\alpha \in R^1$ are, respectively, defined as

$$T_{ijk}(\boldsymbol{\omega}) = \mathrm{E}_{\boldsymbol{\omega}}[\partial_i \ell(\boldsymbol{\omega}|\mathbf{Y}, \boldsymbol{\theta}) \, \partial_j \ell(\boldsymbol{\omega}|\mathbf{Y}, \boldsymbol{\theta}) \, \partial_k \ell(\boldsymbol{\omega}|\mathbf{Y}, \boldsymbol{\theta})],$$

and $\Gamma^\alpha_{ijk}(\boldsymbol{\omega}) = \mathrm{E}_{\boldsymbol{\omega}}[\partial_i \partial_j \ell(\boldsymbol{\omega}|\mathbf{Y}, \boldsymbol{\theta}) \, \partial_k \ell(\boldsymbol{\omega}|\mathbf{Y}, \boldsymbol{\theta})] + 0.5(1-\alpha)T_{ijk}(\boldsymbol{\omega})$. It can be shown that $\Gamma^\alpha_{ijk}(\boldsymbol{\omega}) = \Gamma^0_{ijk}(\boldsymbol{\omega}) - \alpha T_{ijk}(\boldsymbol{\omega})/2$, where $\Gamma^0_{ijk}(\boldsymbol{\omega})$ is the Christoffel symbol for the Lévi–Civita connection of the metric tensor and

$$\Gamma^0_{ijk}(\boldsymbol{\omega}) = \tfrac{1}{2}[\partial_i g_{jk}(\boldsymbol{\omega}) + \partial_j g_{ik}(\boldsymbol{\omega}) - \partial_k g_{ij}(\boldsymbol{\omega})].$$

With the above quantities, the perturbation model $M$ is a statistical manifold, which plays an important role in understanding the behavior of the perturbed model (Amari [1], Kass and Vos [14], Lauritzen [15] and Zhu and Wei [39]).

DEFINITION 2. A *statistical perturbation manifold* $(M, G(\boldsymbol{\omega}), T(\boldsymbol{\omega}))$ is the manifold $M$ with a metric $G(\boldsymbol{\omega})$ and a covariant 3-tensor $T(\boldsymbol{\omega})$.

Now we consider a specific smooth curve in $M$, called an $\alpha$-geodesic.

DEFINITION 3. $\boldsymbol{\omega}(t)$ is called an *$\alpha$-geodesic* with respect to the affine connection $\Gamma^\alpha_{ijk}(\boldsymbol{\omega})$ if it satisfies the equation

$$\frac{d^2\omega_i(t)}{dt^2} + \sum_{s,j,k} g^{is}(\boldsymbol{\omega}(t))\Gamma^\alpha_{jks}(\boldsymbol{\omega}(t)) \frac{d\omega_j(t)}{dt} \frac{d\omega_k(t)}{dt} = 0, \tag{8}$$

where $g^{is}(\boldsymbol{\omega})$ is the $(i,s)$th element of $G(\boldsymbol{\omega})^{-1}$.

The geodesic is a direct extension of the straight line $\boldsymbol{\omega}(t) = \boldsymbol{\omega}^0 + t\mathbf{h}$ in Euclidean space (Amari [1] and Kass and Vos [14]). In particular, as we move along a geodesic, the tangent vector of the geodesic does not change in length and direction. If $\Gamma^\alpha_{ijk}(\boldsymbol{\omega}) = 0$ for all $\boldsymbol{\omega}$, then the manifold is $\alpha$-flat and the geodesic equation for this $\alpha$ is linear in $t$: $\boldsymbol{\omega}(t) = \boldsymbol{\omega}^0 + t\mathbf{h}$.

Some important properties related to the above geometrical quantities are summarized in the following lemma, whose proof can be found in Amari [1], pages 40, 51–52.

LEMMA 1. *Let $\boldsymbol{\phi} = (\phi^1, \ldots, \phi^p) = \boldsymbol{\phi}(\boldsymbol{\omega})$ be a new coordinate system of $M$, $\partial_a = \partial/\partial\phi^a$, $B^a_i = \partial\phi^a/\partial\omega_i$ and $B^i_a = \partial\omega_i/\partial\phi^a$. Then the geometrical quantities of $M$ in the coordinate system $\boldsymbol{\phi}$ can be written as $g_{ab} = \sum_{i,j} B^i_a B^j_b g_{ij}$, $T_{abc} = \sum_{i,j,k} B^i_a B^j_b B^k_c T_{ijk}$ and $\Gamma^\alpha_{abc} = \sum_{i,j,k} B^i_a B^j_b B^k_c \Gamma^\alpha_{ijk} + \sum_{i,j} g_{ij} B^i_c \partial_a B^j_b$. We use the indices $i, j, k$, and so on, to denote quantities related to the coordinate system $\boldsymbol{\omega}$ and the indices $a, b, c$, and so on, to denote quantities related to the coordinate system $\boldsymbol{\phi}$.*

2.3. *Influence measures and their properties.* Let $f(\boldsymbol{\omega})\colon R^p \to R^1$ be the objective function (e.g., the likelihood displacement function in Cook [7] or the residual sum of squares in Wu and Luo [32]), which defines the aspect of inference of interest for sensitivity analysis. Let $\boldsymbol{\omega}(t)$ be a smooth curve on $M$ with $\boldsymbol{\omega}(0) = \boldsymbol{\omega}^0$ and $d\boldsymbol{\omega}(t)/dt|_{t=0} = \mathbf{h} \in T_{\boldsymbol{\omega}^0}$. Therefore, $f(\boldsymbol{\omega}(t))$ is a function of $\boldsymbol{\omega}(t)$ defined on the perturbation manifold $M$. It follows from a Taylor series expansion that

$$(9) \qquad f(\boldsymbol{\omega}(t)) = f(\boldsymbol{\omega}(0)) + \dot{f}_{\mathbf{h}}(0)t + \tfrac{1}{2}\ddot{f}_{\mathbf{h}}(0)t^2 + o(t^2).$$

The first and second derivatives of $f(\boldsymbol{\omega}(t))$ at $t = 0$ are, respectively, given by

$$(10) \quad \dot{f}_{\mathbf{h}}(0) = \sum_j \frac{\partial f(\boldsymbol{\omega}^0)}{\partial\omega_j}h_j = \nabla^T_f \mathbf{h} \quad \text{and} \quad \ddot{f}_{\mathbf{h}}(0) = \mathbf{h}^T H_f \mathbf{h} + \nabla^T_f \frac{d^2\boldsymbol{\omega}(0)}{dt^2}.$$

If $\nabla_f \neq 0$, then the first-order term $\dot{f}_{\mathbf{h}}(0)$ mainly characterizes the local influence of a perturbation vector $\boldsymbol{\omega}$ to a model. However, if $\nabla_f = 0$, then it follows from (9) and (10) that $\dot{f}_{\mathbf{h}}(0)$ and $\ddot{f}_{\mathbf{h}}(0)$ reduce to zero and $\mathbf{h}^T H_f \mathbf{h}$, respectively; therefore, we must use the second order term $\ddot{f}_{\mathbf{h}}(0)$ to assess the local behavior of the objective function when $\nabla_f = 0$.

An important question is how to assess the local influence of minor perturbations to a model when $\nabla_f \neq 0$. This is the so-called first-order approach in Wu and Luo [31]. For instance, in a transformation model, $f$ can be defined as the transformation parameter estimate, whose first derivative does not equal 0 at $\boldsymbol{\omega}^0$ (Lawrance [16]). We introduce a first-order influence measure as follows.

DEFINITION 4. The *first-order influence measure* (FI) in the direction $\mathbf{h} \in T_{\boldsymbol{\omega}^0}$ is defined as

$$\text{(11)} \qquad \text{FI}_{f,\mathbf{h}} = \text{FI}_{f(\boldsymbol{\omega}^0),\mathbf{h}} = \frac{\mathbf{h}^T \nabla_f \nabla_f^T \mathbf{h}}{\mathbf{h}^T G \mathbf{h}},$$

where $G = G(\boldsymbol{\omega}^0)$.

The proposed FI has an interesting geometrical interpretation and is invariant with respect to arbitrary reparametrizations at any point in $\boldsymbol{\omega}$. We are now led to the following theorem.

THEOREM 1. *We have the following results:*

(i) $\text{FI}_{f,\mathbf{h}} = \lim_{t \to 0} [f(\boldsymbol{\omega}(t)) - f(\boldsymbol{\omega}(0))]^2 / S(\boldsymbol{\omega}(0), \boldsymbol{\omega}(t))^2$.

(ii) *If $\boldsymbol{\phi}$ is a diffeomorphism of $\boldsymbol{\omega}$, then $\text{FI}_{f(\boldsymbol{\omega}),\mathbf{h}}$ is invariant with respect to any reparametrization corresponding to $\boldsymbol{\phi}$ and $\text{FI}_{kf,\mathbf{h}} = k^2 \text{FI}_{f,\mathbf{h}}$ holds for any $k$.*

PROOF. It follows from (7) that $S(\boldsymbol{\omega}(0), \boldsymbol{\omega}(t))^2 = t^2 \mathbf{h}^T G \mathbf{h} + o(t^2)$. Using l'Hôpital's rule and (9), we can prove (i). Assuming $\boldsymbol{\omega} = \boldsymbol{\omega}(\boldsymbol{\phi})$ and $\boldsymbol{\phi} = \boldsymbol{\phi}(\boldsymbol{\omega})$, the Jacobian matrices of the above coordinate transformations are given by $\Phi = \partial \boldsymbol{\phi} / \partial \boldsymbol{\omega}$ and $\Psi = \partial \boldsymbol{\omega} / \partial \boldsymbol{\phi}$. Differentiating the identities $\boldsymbol{\phi}[\boldsymbol{\omega}(\boldsymbol{\phi})] = \boldsymbol{\phi}$ and $\boldsymbol{\omega}[\boldsymbol{\phi}(\boldsymbol{\omega})] = \boldsymbol{\omega}$ with respect to $\boldsymbol{\phi}$ and $\boldsymbol{\omega}$, respectively, leads to $\Psi\Phi = \Phi\Psi = \mathbf{I}_p$. Thus, we have $G(\boldsymbol{\phi}) = \Psi^T G(\boldsymbol{\omega})\Psi$ and $\nabla_{f(\boldsymbol{\phi}^0)} = \Psi^T \nabla_{f(\boldsymbol{\omega}^0)}$, where $\boldsymbol{\phi}^0 = \boldsymbol{\omega}^0$. Using Definition 4, we can prove (ii). □

The statistical significance of Theorem 1 is two-fold. First, Theorem 1(i) indicates that the first-order measure is associated with the first derivative of $f(\boldsymbol{\omega}(t))$ on $M$ evaluated at $t = 0$. If $M$ is a Euclidean space, $\mathbf{h}^T \mathbf{h} = 1$ and $\boldsymbol{\omega}(t) = t\mathbf{h} + \boldsymbol{\omega}^0$, then $\text{FI}_{f,\mathbf{h}}$ reduces to the square of the directional derivative of $f$ at $\boldsymbol{\omega}^0$ in the direction $\mathbf{h}$, given by $\lim_{t \to 0} [f(\boldsymbol{\omega}^0 + t\mathbf{h}) - f(\boldsymbol{\omega}^0)]^2 / t^2$. Second, although $\boldsymbol{\omega}$ may not be an appropriate perturbation, we can always use $G$ to obtain an appropriate perturbation $\tilde{\boldsymbol{\omega}}$ in (5), which yields

$$\text{FI}_{f(\tilde{\boldsymbol{\omega}}),\mathbf{h}}|_{\tilde{\boldsymbol{\omega}} = \boldsymbol{\omega}^0} = \frac{\mathbf{h}^T G^{-1/2} \nabla_f \nabla_f^T G^{-1/2} \mathbf{h}}{\mathbf{h}^T \mathbf{h}}.$$

The maximum value of $\text{FI}_{f,\mathbf{h}}$ equals $\nabla_f^T G^{-1} \nabla_f$, which quantifies the degree of local influence of $\tilde{\boldsymbol{\omega}}$ to a statistical model, while the corresponding direction vector $\tilde{\mathbf{h}}_{\max} = G^{-1/2} \nabla_f$ can be used for identifying influential observations (Lawrance [16]).

We use $\ddot{f}_{\mathbf{h}}(0)$ to assess the second-order local influence of $\boldsymbol{\omega}$ to a statistical model, even when $\nabla_f \neq 0$. The approach which utilizes the information in $\ddot{f}_{\mathbf{h}}(0)$ is called the second-order approach (Wu and Luo [31, 32]). However, for a general curve $\boldsymbol{\omega}(t)$ on $M$, $\dot{f}_{\mathbf{h}}(0)$ may not be geometrically well behaved (Murray and Rice [20]). Instead, we only consider the 0-geodesic $\boldsymbol{\omega}(t)$ associated with the Lévi–Civita connection of the metric tensor $G(\boldsymbol{\omega})$, which is unique and defined in an interval containing 0 such that $\boldsymbol{\omega}(t) = \boldsymbol{\omega}^0$ and $d\boldsymbol{\omega}(t)/dt = \mathbf{h} \in T_{\boldsymbol{\omega}^0}$. Then, we can obtain a covariant version of Taylor's theorem (Murray and Rice [20] and McCullagh and Cox [19]) as follows:

$$(12) \qquad f(\boldsymbol{\omega}(t)) = f(\boldsymbol{\omega}^0) + t\nabla_f^T \mathbf{h} + \tfrac{1}{2}t^2 \mathbf{h}^T \tilde{H}_f^0 \mathbf{h} + o(t^2),$$

where $\tilde{H}_f^0 = \tilde{H}_{f(\boldsymbol{\omega}^0)}^0$ and the $(i,j)$th element of $\tilde{H}_{f(\boldsymbol{\omega})}^0$ is given by

$$[\tilde{H}_{f(\boldsymbol{\omega})}^0]_{(i,j)} = \partial_i \partial_j f(\boldsymbol{\omega}) - \sum_{s,r} g^{sr}(\boldsymbol{\omega}) \Gamma_{ijs}^0(\boldsymbol{\omega}) \, \partial_r f(\boldsymbol{\omega}).$$

The matrix $\tilde{H}_{f(\boldsymbol{\omega})}^0$ is called the covariant Hessian of $f(\boldsymbol{\omega})$ (Zhu and Wei [39] and Murray and Rice [20]). Because $\Gamma_{ijk}^0(\boldsymbol{\omega})$ is symmetric with respect to $i$ and $j$, it can be shown that $\tilde{H}_{f(\boldsymbol{\omega})}^0$ is a symmetric matrix. In particular, $\tilde{H}_{f(\boldsymbol{\omega})}^0$ satisfies the following property (Murray and Rice [20]).

LEMMA 2. *Let $\boldsymbol{\phi}$ be a diffeomorphism of $\boldsymbol{\omega}$ with Jacobian matrix $\Psi = \partial \boldsymbol{\omega}/\partial \boldsymbol{\phi}$. Then $\tilde{H}_{f(\boldsymbol{\phi})}^0 = \Psi^T \tilde{H}_{f(\boldsymbol{\omega})}^0 \Psi$.*

Lemma 2 shows that $\tilde{H}_{f(\boldsymbol{\omega})}^0$ is a 2-tensor, so it is geometrically well behaved.

DEFINITION 5. The *second-order influence measure* (SI) in the direction $\mathbf{h} \in T_{\boldsymbol{\omega}^0}$ is defined as

$$(13) \qquad \mathrm{SI}_{f,\mathbf{h}} = \mathrm{SI}_{f(\boldsymbol{\omega}^0),\mathbf{h}} = \frac{\mathbf{h}^T \tilde{H}_f^0 \mathbf{h}}{\mathbf{h}^T G \mathbf{h}}.$$

The *standardized SI* (SSI) in the direction $\mathbf{h} \in T_{\boldsymbol{\omega}^0}$ is defined as

$$(14) \qquad \mathrm{SSI}_{f,\mathbf{h}} = \mathrm{SSI}_{f(\boldsymbol{\omega}^0),\mathbf{h}} = \frac{1}{\|G^{-1}\tilde{H}_f^0\|_M} \frac{\mathbf{h}^T \tilde{H}_f^0 \mathbf{h}}{\mathbf{h}^T G \mathbf{h}}.$$

We can establish the following properties of $\mathrm{SI}_{f,\mathbf{h}}$ and $\mathrm{SSI}_{f,\mathbf{h}}$.

THEOREM 2. *We have the following results:*

(i) $\mathrm{SI}_{f(\boldsymbol{\omega}^0),\mathbf{h}} = \lim_{t \to 0} 2[f(\boldsymbol{\omega}(t)) - f(\boldsymbol{\omega}(0)) - t\nabla_f^T \mathbf{h}]/S(\boldsymbol{\omega}(0), \boldsymbol{\omega}(t))^2.$

(ii) *Suppose $\phi$ is a diffeomorphism of $\boldsymbol{\omega}$. Then* $\mathrm{SI}_{f(\boldsymbol{\omega}^0),\mathbf{h}}$ *and* $\mathrm{SSI}_{f(\boldsymbol{\omega}^0),\mathbf{h}}$ *are invariant with respect to any reparametrization corresponding to $\boldsymbol{\phi}$ at $\boldsymbol{\omega}^0$. Moreover,*

$$\mathrm{SI}_{kf(\boldsymbol{\omega}),\mathbf{h}} = k\mathrm{SI}_{f(\boldsymbol{\omega}),\mathbf{h}} \quad \text{and} \quad \mathrm{SSI}_{kf(\boldsymbol{\omega}),\mathbf{h}} = \mathrm{SSI}_{f(\boldsymbol{\omega}),\mathbf{h}}, \tag{15}$$

*for any $k \neq 0$ and $\boldsymbol{\omega} \in \Omega$.*

(iii) *Let $\{(\lambda_i, \mathbf{u}_i), i = 1, \ldots, p\}$ be the eigenvalue–eigenvector (E–E) pairs of $\tilde{H}_f^0$ with respect to $G$. Then, for any direction $\mathbf{h}$, we have $0 \leq \mathrm{SSI}_{f,\mathbf{h}} \leq 1$,*

$$\mathrm{SI}_{f,\mathbf{u}_i} = \lambda_i \quad \text{and} \quad \mathrm{SSI}_{f,\mathbf{u}_i} = \hat{\lambda}_i = \frac{\lambda_i}{\sqrt{\sum_{j=1}^p \lambda_j^2}},$$

*where $\hat{\lambda}_i$ is the normalized eigenvalue.*

PROOF. Using (7), (12) and l'Hôpital's rule, we can prove Theorem 2(i). Because a diffeomorphism exists between $\boldsymbol{\omega}$ and $\boldsymbol{\phi}$ such that $\boldsymbol{\omega} = \boldsymbol{\omega}(\boldsymbol{\phi})$ and $\boldsymbol{\phi} = \boldsymbol{\phi}(\boldsymbol{\omega})$, we have $\Psi\Phi = \Phi\Psi = \mathbf{I}_p$. Moreover, because $G(\boldsymbol{\phi})$ is a metric tensor and $\tilde{H}_{f(\boldsymbol{\phi})}^0$ is a 2-tensor, we have

$$G(\boldsymbol{\phi}) = \Psi^T G(\boldsymbol{\omega})\Psi \quad \text{and} \quad \tilde{H}_{f(\boldsymbol{\phi}^0)}^0 = \Psi^T \tilde{H}_{f(\boldsymbol{\omega}^0)}^0 \Psi.$$

Consider a geodesic $\boldsymbol{\omega}(t)$ with $\boldsymbol{\omega}(0) = \boldsymbol{\omega}^0$ and $d\boldsymbol{\omega}(0)/dt = \mathbf{h} \in T_{\boldsymbol{\omega}^0}$. Then $\boldsymbol{\phi}(\boldsymbol{\omega}(t))$ is a geodesic in the $\boldsymbol{\phi}$-coordinate such that $\boldsymbol{\phi}(\boldsymbol{\omega}^0) = \boldsymbol{\phi}^0$ and $d\boldsymbol{\phi}(\boldsymbol{\omega}(0))/dt = \Phi\mathbf{h}$. If $G$ is positive definite, then it follows from Lemmas 1 and 2 that

$$\mathrm{SI}_{f(\boldsymbol{\phi}^0),\Phi\mathbf{h}} = \frac{\mathbf{h}^T \Phi^T \tilde{H}_{f(\boldsymbol{\phi}^0)}^0 \Phi\mathbf{h}}{\mathbf{h}^T \Phi^T G(\boldsymbol{\phi}^0)\Phi\mathbf{h}} = \frac{\mathbf{h}^T \Phi^T \Psi^T \tilde{H}_{f(\boldsymbol{\omega}^0)}^0 \Psi\Phi\mathbf{h}}{\mathbf{h}^T \Phi^T \Psi^T G(\boldsymbol{\omega}^0)\Psi\Phi\mathbf{h}} = \mathrm{SI}_{f(\boldsymbol{\omega}^0),\mathbf{h}}.$$

Similarly, we can show $\mathrm{SSI}_{f(\boldsymbol{\phi}^0),\Phi\mathbf{h}} = \mathrm{SSI}_{f(\boldsymbol{\omega}^0),\mathbf{h}}$. Thus, $\mathrm{SI}_{f,\mathbf{h}}$ and $\mathrm{SSI}_{f,\mathbf{h}}$ are invariant with respect to reparametrization $\boldsymbol{\phi}$ at $\boldsymbol{\omega}^0$. For any $k$, we have $\tilde{H}_{kf(\boldsymbol{\omega})}^0 = k\tilde{H}_{f(\boldsymbol{\omega})}^0$ and equation (15) holds for any $\boldsymbol{\omega}$ and $k \neq 0$. This proves Theorem 2(ii). By using Definition 5, we thus prove Theorem 2(iii). This completes the proof of Theorem 2. □

Theorem 2 has the following implications. First, if $\boldsymbol{\omega}$ is an appropriate perturbation and $\nabla_f = \mathbf{0}$, then $\mathrm{SSI}_{f,\mathbf{h}} = B_\mathbf{h}$ and $\mathrm{SI}_{f,\mathbf{h}} = C_\mathbf{h}$. We note that most of the examples in Cook [7] fall into this scenario. In general, even though we may choose a perturbation $\boldsymbol{\omega}$ which is not appropriate, we can always use $G$ to obtain an appropriate perturbation $\tilde{\boldsymbol{\omega}}$ in (5). In this case, the normal curvature and the second-order influence measures will lead to the same results when $\nabla_f = \mathbf{0}$ and the chosen perturbation is appropriate. Therefore, the diagnostic method proposed here can be regarded as an extension of Cook's [7] local influence approach in a more general setting. Second, $\mathrm{SI}_{f(\boldsymbol{\omega}),\mathbf{h}}$ and $\mathrm{SSI}_{f(\boldsymbol{\omega}),\mathbf{h}}$ are scale invariant even when $\nabla_f \neq \mathbf{0}$, whereas $C_\mathbf{h}$ and $B_\mathbf{h}$ are not (Fung and Kwan [11]). This generalization facilitates new methods and techniques for doing sensitivity analyses of a statistical model.

2.4. *New local influence approach.* What follows are the four key steps in assessing local influence of perturbing a parametric model $p(\mathbf{Y}|\boldsymbol{\theta})$:

Step 1. Choose a perturbation scheme $\boldsymbol{\omega}$ such that $\int p(\mathbf{Y}|\boldsymbol{\theta},\boldsymbol{\omega})\,d\mathbf{Y} = 1$.

Step 2. Given the perturbed model, calculate the geometrical quantities [e.g., $g_{ij}(\boldsymbol{\omega})$, $T_{ijk}(\boldsymbol{\omega})$, and $\Gamma^{\alpha}_{ijk}(\boldsymbol{\omega})$] of the perturbation manifold.

Step 3. Check whether the perturbation $\boldsymbol{\omega}$ is appropriate, that is, $G(\boldsymbol{\omega}^0) = c\mathbf{I}_p$. If yes, go to Step 4 below. Otherwise, find a new perturbation scheme and go back to Step 2.

Step 4. Choose an objective function $f(\boldsymbol{\omega})$. If $\nabla_f = \mathbf{0}$, then use SI and SSI to assess local influence of minor perturbations to a model. However, if $\nabla_f$ is nonzero, then use FI, SI and SSI together.

**3. Appropriate perturbations in four examples.** We examine four examples to illustrate how to calculate geometrical quantities for a perturbation manifold and show how to find an appropriate perturbation in the examples. We also consider several objective functions to assess the local influence of an appropriate perturbation to a parametric model in each of the four examples.

3.1. *Case-weight perturbation.* Suppose that $Y_1, \ldots, Y_n$ are independent with $m_1 = \cdots = m_n = 1$ and $p(\mathbf{y}|\boldsymbol{\theta})$ can be written as $\prod_{i=1}^n p(y_i; \boldsymbol{\theta})$. We consider case-weight perturbation, in which $L(\boldsymbol{\theta}|\boldsymbol{\omega})$ is given by

$$(16) \qquad L(\boldsymbol{\theta}|\boldsymbol{\omega}) = \sum_{i=1}^n \omega_i \ell(y_i; \boldsymbol{\theta}),$$

where $\ell(y_i; \boldsymbol{\theta}) = \log p(y_i; \boldsymbol{\theta})$. Thus, $p = n$ and $\boldsymbol{\omega}^0 = \mathbf{1}_n$ is an $n \times 1$ vector with all elements equal to 1. The density of the perturbation model $p(\mathbf{y}|\boldsymbol{\theta},\boldsymbol{\omega})$ is given by

$$(17) \qquad p(\mathbf{y}|\boldsymbol{\theta},\boldsymbol{\omega}) = \prod_{i=1}^n \{\exp\{\omega_i \ell(y_i; \boldsymbol{\theta})\}[c_i(\omega_i; \boldsymbol{\theta})]^{-1}\},$$

where $c_i(\omega_i; \boldsymbol{\theta}) = \int \exp\{\omega_i \ell(y_i; \boldsymbol{\theta})\}\,dy_i$ for all $i = 1, \ldots, n$. After some algebraic calculations, we have the following results.

THEOREM 3. *If the four regularity conditions of Amari [1] hold for $p(\mathbf{y}|\boldsymbol{\theta},\boldsymbol{\omega})$, then the following results hold for case-weight perturbation:*

(i) $\partial_i \log c_i(\omega_i; \boldsymbol{\theta}) = \mathrm{E}_\omega[\ell(y_i; \boldsymbol{\theta})]$ *and* $\partial_\theta \log c_i(\omega_i; \boldsymbol{\theta}) = \omega_i \mathrm{E}_\omega[\partial_\theta \ell(y_i; \boldsymbol{\theta})]$, *where the expectation* $\mathrm{E}_\omega$ *is taken with respect to* $p(\mathbf{y}|\boldsymbol{\theta},\boldsymbol{\omega})$;

(ii) $\partial_i^2 \log c_i(\omega_i; \boldsymbol{\theta}) = \mathrm{Var}_\omega[\ell(y_i; \boldsymbol{\theta})]$,

$\qquad \partial_\theta^2 \log c_i(\omega_i; \boldsymbol{\theta}) = \omega_i \mathrm{E}_\omega[\partial_\theta^2 \ell(y_i; \boldsymbol{\theta})] + \omega_i^2 \mathrm{Var}_\omega[\partial_\theta \ell(y_i; \boldsymbol{\theta})]$

*and*

$$\partial_\theta \partial_i \log c_i(\omega_i; \boldsymbol{\theta})$$
$$= \mathrm{E}_\omega[\partial_\theta \ell(y_i; \boldsymbol{\theta})]$$
$$+ \omega_i \mathrm{E}_\omega(\{\partial_\theta \ell(y_i; \boldsymbol{\theta}) - \mathrm{E}_\omega[\partial_\theta \ell(y_i; \boldsymbol{\theta})]\}\{\ell(y_i; \boldsymbol{\theta}) - \mathrm{E}_\omega[\ell(y_i; \boldsymbol{\theta})]\});$$

(iii) $g_{ij}(\boldsymbol{\omega}) = \mathrm{Var}_\omega[\ell(y_i; \boldsymbol{\theta})]\delta_{ij}$, $\Gamma_{ijk}^\alpha(\boldsymbol{\omega}) = 0.5(1-\alpha)T_{ijk}(\boldsymbol{\omega})$ *and*

$$T_{ijk}(\boldsymbol{\omega}) = \mathrm{E}_\omega[\{\ell(y_i; \boldsymbol{\theta}) - \mathrm{E}_\omega[\ell(y_i; \boldsymbol{\theta})]\}^3]\delta_{ij}\delta_{ik} \qquad for \ i,j,k = 1, \ldots, n,$$

*where $\delta_{ij}$ is the Kronecker delta;*

(iv) *the geodesic $\boldsymbol{\omega} = \boldsymbol{\omega}(t)$ with respect to $\Gamma_{ijk}^\alpha(\boldsymbol{\omega})$ satisfies*

$$\int_{\boldsymbol{\omega}^0}^{\boldsymbol{\omega}} \exp\left\{\int g^{ii}(\boldsymbol{\xi})\Gamma_{iii}^\alpha(\boldsymbol{\xi})\, d\boldsymbol{\xi}\right\} d\boldsymbol{\xi} = h_i t$$

*for $i = 1, \ldots, n$, where $\boldsymbol{\omega}(0) = \boldsymbol{\omega}^0$ and $d\boldsymbol{\omega}(t)/dt = \mathbf{h} = (h_1, \ldots, h_n)^T$. In particular, $\boldsymbol{\omega}(t) = \boldsymbol{\omega}^0 + t\mathbf{h}$ is a 1-geodesic.*

PROOF. The log-likelihood function of $p(\mathbf{y}|\boldsymbol{\theta}, \boldsymbol{\omega})$ is given by

$$\ell(\boldsymbol{\omega}|\mathbf{y}, \boldsymbol{\theta}) = \log p(\mathbf{y}|\boldsymbol{\eta}) = \sum_{i=1}^n [\omega_i \ell(y_i; \boldsymbol{\theta}) - \log c_i(\omega_i; \boldsymbol{\theta})],$$

where $\boldsymbol{\eta} = (\boldsymbol{\theta}, \boldsymbol{\omega})$. By using $\mathrm{E}_\omega[\partial_\eta \ell(\boldsymbol{\omega}|\mathbf{y}, \boldsymbol{\theta})] = 0$ and $\mathrm{E}_\omega[-\partial_\eta^2 \ell(\boldsymbol{\omega}|\mathbf{y}, \boldsymbol{\theta})] = \mathrm{E}_\omega\{[\partial_\eta \ell(\boldsymbol{\omega}|\mathbf{y}, \boldsymbol{\theta})][\partial_\eta \ell(\boldsymbol{\omega}|\mathbf{y}, \boldsymbol{\theta})]^T\}$, we can obtain (i) and (ii).

By differentiating $\ell(\boldsymbol{\omega}|\mathbf{y}, \boldsymbol{\theta})$ with respect to $\boldsymbol{\omega}$, we have $\partial_i \ell(\boldsymbol{\omega}|\mathbf{y}, \boldsymbol{\theta}) = \ell(y_i; \boldsymbol{\theta}) - \partial_i \log c_i(\omega_i; \boldsymbol{\theta})$ and $\partial_i^2 \ell(\boldsymbol{\omega}|\mathbf{y}, \boldsymbol{\theta}) = -\partial_i^2 \log c_i(\omega_i; \boldsymbol{\theta})$. Therefore, we can directly calculate the geometric quantities $g_{ij}(\boldsymbol{\omega})$, $\Gamma_{ijk}^\alpha(\boldsymbol{\omega})$ and $T_{ijk}(\boldsymbol{\omega})$, which lead to (iii).

The geodesic $\boldsymbol{\omega} = \boldsymbol{\omega}(t)$ with respect to $\Gamma_{ijk}^\alpha(\boldsymbol{\omega})$ satisfies a second-order differential equation which is defined by $d^2\omega_i(t)/dt^2 + g^{ii}(\boldsymbol{\omega})\Gamma_{iii}^0(\boldsymbol{\omega})[d\omega_i(t)/dt]^2 = 0$, with initial conditions $\boldsymbol{\omega}(0) = \boldsymbol{\omega}^0$ and $d\boldsymbol{\omega}(t)/dt = \mathbf{h}$. We can prove (iv) by solving this second-order differential equation (Coddington [6]).

Theorem 3 establishes the manifold of case-weight perturbation in (16). If the $Y_i$ are also identically distributed, then $G(\boldsymbol{\omega}^0) = g_{11}(\boldsymbol{\omega}^0)\mathbf{I}_n$ and the perturbation in (16) is an appropriate one. Moreover, if we treat both responses and covariates as random variables, the perturbation in (16) is appropriate even for the regression case. In general, if

$$G = \mathrm{diag}(\mathrm{Var}_{\boldsymbol{\omega}^0}[\ell(y_1; \boldsymbol{\theta})], \ldots, \mathrm{Var}_{\boldsymbol{\omega}^0}[\ell(y_n; \boldsymbol{\theta})]) \neq c\mathbf{I}_n$$

for any $c > 0$, then we consider a new perturbation $\tilde{\boldsymbol{\omega}}$ in (5) with $c = 1$ such that $G(\tilde{\boldsymbol{\omega}}) = \mathbf{I}_n$ at $\boldsymbol{\omega}^0$ and

(18)    $$L(\boldsymbol{\theta}|\tilde{\boldsymbol{\omega}}) = \sum_{i=1}^n [1 + \sqrt{\mathrm{Var}_{\boldsymbol{\omega}^0}[\ell(y_i; \boldsymbol{\theta})]}(\tilde{\omega}_i - 1)]\ell(y_i; \boldsymbol{\theta}).$$

With the development above, we can now choose an objective function (e.g., the likelihood displacement) and calculate its associated influence measures {FI, SI, SSI} to assess local influence of the perturbation (16). For instance, if we are interested in a particular component of $\theta$, say, $\theta_1$, we may use $\hat{\theta}_1(\boldsymbol{\omega})$ as an objective function, where $\hat{\theta}_1(\boldsymbol{\omega})$ is the maximum likelihood estimate of $\theta_1$ under the perturbation (18). $\square$

3.2. *Location-scale family.* Suppose that $Y_1, \ldots, Y_n$ are independent and $m_1 = \cdots = m_n = 1$. Let $\boldsymbol{\theta} = (\boldsymbol{\beta}^T, \sigma^2)$. Each $p(y_i; \boldsymbol{\theta}) = \sigma^{-1} p_0(\sigma^{-1}(y_i - \mu(\mathbf{x}_i^T \boldsymbol{\beta})))$ belongs to a location-scale family such that $p_0$ is a known density satisfying $\int x p_0(x) = 0$ and $\int x^2 p_0(x) = 1$, where $\mu(\cdot)$ is a given function and $\mathbf{x}_i$ is a $q_1 \times 1$ vector. Thus, $\mathrm{E}(y_i) = \mu_i = \mu(\mathbf{x}_i^T \boldsymbol{\beta})$ and $\mathrm{Var}(y_i) = \sigma^2$.

We consider three different perturbations: case-weight perturbation, perturbation of the variance and perturbation of the response. Then we establish a perturbation manifold for each of these perturbations as follows.

For case-weight perturbation, since $\ell(y_i; \boldsymbol{\theta}) = \log p_0(\sigma^{-1}(y_i - \mu_i)) - \log \sigma$, we can use the transformation $e_i = (y_i - \mu_i)/\sigma$ to obtain $g_{ij}(\boldsymbol{\omega}) = \mathrm{Var}_{\omega,0}[\log p_0(e_i)]\delta_{ij}$. Similarly, we can calculate $\Gamma_{ijk}^\alpha(\boldsymbol{\omega}) = (1 - \alpha)T_{ijk}(\boldsymbol{\omega})$ and $T_{ijk}(\boldsymbol{\omega}) = \mathrm{E}_{\omega,0}(\{\log p_0(e_i) - \mathrm{E}_{\omega,0}[\log p_0(e_i)]\}^2)\delta_{ij}\delta_{kj}$.

For the perturbation of variance, we consider a heterogeneous variance of $\mathbf{y} = (y_1, \ldots, y_n)^T$ such that $\mathrm{Var}(\mathbf{y}) = \sigma^2 \mathrm{diag}(1/\omega_1^2, \ldots, 1/\omega_n^2)$. In this case, $\boldsymbol{\omega}^0 = \mathbf{1}_n$ and $p = n$. The log-likelihood function for the perturbed model is given by $-n \log \sigma + \sum_{i=1}^n \log \omega_i + \sum_{i=1}^n \log p_0(t_i)$, where $t_i = \omega_i \sigma^{-1}[y_i - \mu(\mathbf{x}_i^T \boldsymbol{\beta})]$. After some algebraic manipulations, we get $g_{ij}(\boldsymbol{\omega}) = \delta_{ij}\omega_i^{-2}\mathrm{E}_0\{[x\partial_x \log p_0(x) + 1]^2\}$, $T_{ijk}(\boldsymbol{\omega}) = \delta_{ij}\delta_{ik}\omega_i^{-3}\mathrm{E}_0\{[x\partial_x \log p_0(x) + 1]^3\}$, and

$$\Gamma_{ijk}^\alpha(\boldsymbol{\omega}) = -\delta_{ij}\delta_{ik}\omega_i^{-3}(\mathrm{E}_0\{[x\partial_x \log p_0(x) + 1]^2\}$$
$$+ 0.5\alpha \mathrm{E}_0\{[x\partial_x \log p_0(x) + 1]^3\}),$$

where $\partial_x = \partial/\partial x$ and the expectation $\mathrm{E}_0$ is taken with respect to $p_0(x)$.

For the perturbation of response, we consider adding a perturbation $\omega_i$ to $y_i$ such that $\ell(\boldsymbol{\omega}|\mathbf{Y}, \boldsymbol{\theta})$ is given by $-n \log \sigma + \sum_{i=1}^n \log p_0(\sigma^{-1}[y_i + \omega_i - \mu(\mathbf{x}_i^T \boldsymbol{\beta})])$. Let $t_i = \sigma^{-1}[y_i + \omega_i - \mu(\mathbf{x}_i^T \boldsymbol{\beta})]$. With some calculations, we can obtain $g_{ij}(\boldsymbol{\omega}) = \delta_{ij}\sigma^{-2}\mathrm{E}_0\{[\partial_x \log p_0(x)]^2\}$, $T_{ijk}(\boldsymbol{\omega}) = \delta_{ij}\delta_{ik}\sigma^{-3}\mathrm{E}_0\{[x\partial_x \log p_0(x) + 1]^3\}$ and $\Gamma_{ijk}^0(\boldsymbol{\omega}) = -0.5\alpha T_{ijk}(\boldsymbol{\omega})$. Thus, since $\Gamma_{ijk}^0(\boldsymbol{\omega})$ vanishes for all $i, j, k = 1, \ldots, n$, the straight line $\boldsymbol{\omega}(t) = \boldsymbol{\omega}^0 + t\mathbf{h}$ is a 0-geodesic.

Combining the above results, we have the following theorem.

THEOREM 4. *If all $p(y_i; \boldsymbol{\theta}) = \sigma^{-1} p_0(\sigma^{-1}(y_i - \mu(\mathbf{x}_i^T \boldsymbol{\beta})))$ belong to the same location-scale family, then*

$$G(\boldsymbol{\omega}^0) = c\mathbf{I}_n \quad and \quad C_{\mathbf{h}} = c\mathrm{SI}_{LD(\boldsymbol{\omega}^0), \mathbf{h}}$$

*hold for case-weight perturbation, the perturbation of variance and the perturbation of response, where $LD(\boldsymbol{\omega})$ denotes the likelihood displacement function in Cook [7].*

Theorem 4 indicates that the three perturbation schemes considered here are appropriate perturbations. Therefore, for $LD(\boldsymbol{\omega})$ introduced in Cook [7], both {SI, SSI} and the normal curvature lead to the same results for location-scale families under the three commonly used perturbations.

### 3.3. *Linear regression model.* Consider the linear regression model

$$(19) \qquad\qquad Y = X\boldsymbol{\beta} + \boldsymbol{\epsilon},$$

where $Y$ is an $n \times 1$ vector of responses, $X$ is an $n \times q_1$ covariate matrix, $\boldsymbol{\beta}$ is a $q_1 \times 1$ vector of unknown parameters and $\boldsymbol{\epsilon} = (\epsilon_1, \ldots, \epsilon_n)^T$ is an $n \times 1$ random vector of errors with distribution $N[\mathbf{0}, \sigma^2 \mathbf{I}_n]$. Here $\boldsymbol{\theta} = (\boldsymbol{\beta}, \sigma^2)$. In this section we fix $\boldsymbol{\theta}$ at its maximum likelihood estimate.

### 3.3.1. *Perturbation of error variances.* We consider a perturbation to error variances via an $n \times 1$ perturbation vector $\boldsymbol{\omega}$ such that

$$(20) \qquad\qquad \mathrm{Var}(\boldsymbol{\epsilon}) = \sigma^2 \, \mathrm{diag}(1/\omega_1, \ldots, 1/\omega_n)$$

and $\boldsymbol{\omega}^0 = \mathbf{1}_n$. It has been shown in Theorem 4 that $g_{ij}(\boldsymbol{\omega}) = 0.5\omega_i^{-2}\delta_{ij}$, $T_{ijk}(\boldsymbol{\omega}) = -\omega_i^{-3}\delta_{ij}\delta_{jk}$ and $\Gamma_{ijk}^{\alpha}(\boldsymbol{\omega}) = 0.5(1-\alpha)T_{ijk}(\boldsymbol{\omega})$. Thus, the perturbation (20) is an appropriate one. However, for illustrative purposes, assume that we consider another perturbation scheme $\boldsymbol{\phi} = (\phi_1, \ldots, \phi_n)^T$ such that

$$(21) \quad \mathrm{Var}(\epsilon_1) = \sigma^2 \frac{k_0}{k_0 - 1 + \phi_1}, \qquad \mathrm{Var}(\epsilon_2) = \sigma^2 \phi_2^{-1}, \ldots, \mathrm{Var}(\epsilon_n) = \sigma^2 \phi_n^{-1},$$

where $k_0 > 0$. Therefore, it can be shown that $\boldsymbol{\phi}^0 = \boldsymbol{\omega}^0$ and $G(\boldsymbol{\phi}^0) = \mathrm{diag}(1/2k_0^2, 1/2, \ldots, 1/2)$. Thus, $G(\boldsymbol{\phi}^0) = c\mathbf{I}_n$ if and only if $k_0 = 1$. That is, the perturbation vector $\boldsymbol{\phi}$ is appropriate only when $k_0 = 1$ in (21), which reduces to (20).

The perturbation of error variances is applied to the residual sum of squares, that is, $f(\boldsymbol{\omega}) = -RSS(\boldsymbol{\omega}) = -r(\boldsymbol{\omega})^T r(\boldsymbol{\omega})$, where $r(\boldsymbol{\omega})$ is the residual vector under the perturbation (20). It can be shown that at $\boldsymbol{\omega}^0 = \mathbf{1}_n$, $\nabla_{-RSS} = (-r_1^2, \ldots, -r_n^2)^T$ and $H_{-RSS} = 2D(\mathbf{r})P_X D(\mathbf{r})$, where $P_X = (p_{ij}) = X(X^T X)^{-1}X^T$ and $D(\mathbf{r}) = \mathrm{diag}(r_1, \ldots, r_n)$, in which $r_i, i = 1, \ldots, n$, are ordinary residuals when $\boldsymbol{\omega}^0 = \mathbf{1}_n$. Because $G = 0.5\mathbf{I}_n$, $\mathrm{FI}_{-RSS,\mathbf{h}} = 2(\mathbf{h}^T \nabla_{-RSS})^2/\mathbf{h}^T\mathbf{h}$ for any $\mathbf{h}$. In particular, we can obtain the maximum value of $\mathrm{FI}_{-RSS,\mathbf{h}}$ as $\sum_{i=1}^{n} r_i^4$ and the corresponding direction vector is given by $\mathbf{h}_{\max} = \nabla_{-RSS}/\sqrt{\sum_{j=1}^{n} r_j^4}$, which is the same as Lawrance's [16] diagnostic. Subsequently, we use $\mathrm{SI}_{-RSS,\mathbf{h}}$ to assess the local influence of the perturbation (20) to model (19). For $\alpha = 0$, it can be shown that $\tilde{H}^0_{-RSS} = 2D(\mathbf{r})P_X D(\mathbf{r}) - D(\mathbf{r}^2)$, where $D(\mathbf{r}^2) = \text{diag}(r_1^2, \ldots, r_n^2)$. Thus, we have $\text{SI}_{-RSS, \mathbf{E}_i} = 2r_i^2(2p_{ii} - 1.0)\sigma^{-2} \approx -2r_i^2\sigma^{-2}$ for $i = 1, \ldots, n$, where $\mathbf{E}_i$ is an $n \times 1$ vector with $i$th element one and zero otherwise for $i = 1, \ldots, n$. Therefore, if an observation has a large absolute residual, then it will be identified as influential by using $\mathbf{h}_{\max}$ and $\text{SI}_{-RSS, \mathbf{E}_i}$.

### 3.3.2. *Perturbation of the explanatory vector.* Consider the perturbation

$$X_\omega = X + WS, \tag{22}$$

where $W = (\omega_{ik})$ is an $n \times q_1$ matrix of perturbations, $S = \text{diag}(s_1, \ldots, s_{q_1})$ with $s_i \neq 0$ for each $i$ and $s_k$ converts $\omega_{ik}$ to the appropriate size and unit so that $s_k\omega_{ik}$ is compatible with the $(i, k)$th element of $X$. With some calculations, it can be shown that $G(\boldsymbol{\omega}) = S\boldsymbol{\beta}\boldsymbol{\beta}^T S \otimes \mathbf{I}_n/\sigma^2$, where $\otimes$ is the Kronecker product of two matrices. Because $\partial\ell(\boldsymbol{\omega}|\mathbf{Y}, \boldsymbol{\theta})/\partial\omega_{ik}$ for $k = 1, \ldots, q_1$ are linearly dependent, the metric matrix is only positive semidefinite. In this case, the singularity of $G(\boldsymbol{\omega})$ indicates that we have introduced too many perturbation parameters; and therefore, some perturbation parameters should be removed. One possibility is to transform $\boldsymbol{\omega} \in R^p$ into $\boldsymbol{\omega}^* \in R^r$ ($p > r$) and make the metric tensor positive definite in the new coordinate system.

Another perturbation scheme is

$$X_\omega = X + W^*\mathbf{1}^*S, \tag{23}$$

where $W^* = \text{diag}(\omega_1, \ldots, \omega_n)$ and $\mathbf{1}^*$ is an $n \times q_1$ matrix with all elements equal to 1. After some calculations, we can obtain $g_{ij}(\boldsymbol{\omega}) = \delta_{ij}(\sum_{k=1}^{q_1} s_k\beta_k)^2/\sigma^2$ and $T_{ijk}(\boldsymbol{\omega}) = \Gamma^\alpha_{ijk}(\boldsymbol{\omega}) = 0$ for all $i, j, k = 1, \ldots, n$. The corresponding perturbation manifold is $\alpha$-flat for any $\alpha \in R^1$ (Amari [1]). Moreover, the commonly considered line $\boldsymbol{\omega}(t) = \boldsymbol{\omega}^0 + t\mathbf{h}$ is a geodesic with respect to $\Gamma^\alpha_{ijk}(\boldsymbol{\omega})$ for any $\alpha \in R^1$ in this perturbation manifold. In particular, the perturbation (23) is an appropriate one, because the metric matrix $G(\boldsymbol{\omega}) = c\mathbf{I}_n$ and is independent of $\boldsymbol{\omega}$.

### 3.4. *Testing a parametric family.* Suppose that $Y_1, \ldots, Y_n$ are independent and identically distributed and the density of $Y_i$, denoted by $p_0(\mathbf{y}_i; \boldsymbol{\theta})$, belongs to a certain parametric family, such as a Gaussian distribution. We consider a log-linear expansion perturbation of $p_0(\mathbf{y}_i; \boldsymbol{\theta})$ as follows (Claeskens and Hjort [5]). Let $\{\psi_j(\mathbf{y}_i; \boldsymbol{\theta}): j = 1, \ldots, m\}$ be a set of functions such that they are orthogonal with respect to $p_0(\mathbf{y}_i; \boldsymbol{\theta})$ and also orthogonal to $\psi_0(\mathbf{y}_i; \boldsymbol{\theta}) = 1$. That is, $\int \psi_j(\mathbf{y}_i; \boldsymbol{\theta})\psi_k(\mathbf{y}_i; \boldsymbol{\theta})p_0(\mathbf{y}_i; \boldsymbol{\theta}) \, d\mathbf{y}_i = \delta_{jk}\text{E}_0(\psi_k^2)$ for all $j, k = 0, \ldots, m$

and $i = 1, \ldots, n$, where $\mathrm{E}_0$ denotes the expectation with respect to $p_0(\mathbf{y}_i; \boldsymbol{\theta})$. For each $i$, the log-linear expansion perturbation is defined by

$$(24) \qquad p(\mathbf{y}_i; \boldsymbol{\theta}, \boldsymbol{\omega}) = p_0(\mathbf{y}_i; \boldsymbol{\theta}) c(\boldsymbol{\omega})^{-1} \exp\left\{ \sum_{j=1}^{m} \omega_j \psi_j(\mathbf{y}_i; \boldsymbol{\theta}) \right\},$$

where $c(\boldsymbol{\omega}) = \int p_0(\mathbf{y}_i; \boldsymbol{\theta}) \exp\{\sum_{j=1}^{m} \omega_j \psi_j(\mathbf{y}_i; \boldsymbol{\theta})\} \, d\mathbf{y}_i$. Thus, the log-likelihood function for the perturbed model is given by

$$\ell(\boldsymbol{\omega}|\mathbf{Y}, \boldsymbol{\theta}) = \sum_{i=1}^{n} \log p_0(\mathbf{y}_i; \boldsymbol{\theta}) + n \sum_{j=1}^{m} \omega_j \overline{\psi_j}(\boldsymbol{\theta}) - n\phi(\boldsymbol{\omega}),$$

where $\overline{\psi_j}(\boldsymbol{\theta}) = \sum_{i=1}^{n} \psi_j(\mathbf{y}_i; \boldsymbol{\theta})/n$ and $\phi(\boldsymbol{\omega}) = \log c(\boldsymbol{\omega})$. In this case, $\boldsymbol{\omega}^0 = \mathbf{0}_m$, an $m \times 1$ vector with all elements 0, and $p = m$.

After some algebraic derivations, we can obtain the geometrical quantities of the perturbation (24) as follows. Because $\partial_j \ell(\boldsymbol{\omega}|\mathbf{Y}, \boldsymbol{\theta}) = n\overline{\psi_j}(\boldsymbol{\theta}) - n \partial_j \phi(\boldsymbol{\omega})$, we have $g_{jk}(\boldsymbol{\omega}) = n \partial_j \partial_k \phi(\boldsymbol{\omega})$, $T_{jkl}(\boldsymbol{\omega}) = n \partial_j \partial_k \partial_l \phi(\boldsymbol{\omega})$ and $\Gamma_{jkl}^{\alpha}(\boldsymbol{\omega}) = 0.5(1-\alpha)T_{jkl}(\boldsymbol{\omega})$ (Amari [1]). In particular, $G(\boldsymbol{\omega}^0) = n \operatorname{diag}(\mathrm{E}_0(\psi_1^2), \ldots, \mathrm{E}_0(\psi_m^2))$. It can be shown that the perturbation in (24) is an appropriate one if and only if $\mathrm{E}_0(\psi_j^2)$ are homogenous, that is, $\mathrm{E}_0(\psi_1^2) = \cdots = \mathrm{E}_0(\psi_m^2)$. Even though $\mathrm{E}_0(\psi_j^2)$ are not homogenous, we can always choose a new perturbation $\tilde{\boldsymbol{\omega}} = G(\boldsymbol{\omega}^0)^{1/2}\boldsymbol{\omega}$ such that

$$(25) \qquad p(\mathbf{y}_i; \boldsymbol{\theta}, \boldsymbol{\omega}) = p_0(\mathbf{y}_i; \boldsymbol{\theta}) c(\tilde{\boldsymbol{\omega}})^{-1} \exp\left\{ \sum_{j=1}^{m} \tilde{\omega}_j \psi_j(\mathbf{y}_i; \boldsymbol{\theta}) / \sqrt{\mathrm{E}_0(\psi_j^2)} \right\}.$$

This $\tilde{\boldsymbol{\omega}}$ is an appropriate perturbation with $G(\tilde{\boldsymbol{\omega}})|_{\mathbf{0}_m} = \mathbf{I}_m$.

We consider the log-likelihood ratio $f(\tilde{\boldsymbol{\omega}}) = \ell(\tilde{\boldsymbol{\omega}}|\mathbf{Y}, \boldsymbol{\theta}) - \log p(\mathbf{Y}|\boldsymbol{\theta})$ as our objective function. Direct calculation leads to

$$\nabla_f = \sqrt{n}\left( \frac{\overline{\psi_1}(\boldsymbol{\theta})}{\sqrt{\mathrm{E}_0(\psi_1^2)}}, \ldots, \frac{\overline{\psi_m}(\boldsymbol{\theta})}{\sqrt{\mathrm{E}_0(\psi_m^2)}} \right)^T \quad \text{and} \quad \mathrm{FI}_{f,\mathbf{h}} = \frac{\mathbf{h}^T \nabla_f \nabla_f \mathbf{h}}{\mathbf{h}^T \mathbf{h}}.$$

The maximum value of $\mathrm{FI}_{f,\mathbf{h}}$ is $\nabla_f^T \nabla_f = n \sum_{j=1}^{m} \overline{\psi_j}(\theta)^2 / \mathrm{E}_0(\psi_j^2)$, which is the well-known score test statistic for testing $H_0: \boldsymbol{\omega} = \mathbf{0}$ when $\theta$ is either known or replaced by its estimate (Claeskens and Hjort [5]). Moreover, for each $\mathbf{E}_j$, $\mathrm{FI}_{f,\mathbf{E}_j} = n\overline{\psi_j}(\theta)^2 / \mathrm{E}_0(\psi_j^2)$ is the score test statistic for testing $H_0: \omega_j = 0$, where $\mathbf{E}_j$ is an $m \times 1$ vector with $i$th element equal to one and zero otherwise for $j = 1, \ldots, m$. We can use $\mathrm{FI}_{f,\mathbf{E}_j}$ to detect the most influential perturbation from all $m$ perturbations. Under some conditions (e.g., $m$ grows slowly with $n$), the score statistic $\nabla_f^T \nabla_f$ converges weakly to a nondegenerate random variable (Claeskens and Hjort [5]), which can be used to characterize the asymptotic behavior of the influence measures $\mathrm{FI}_{f,\mathbf{h}}$.

Combining the above results, we have the following theorem.

THEOREM 5. *If $p_0(\mathbf{y}; \boldsymbol{\theta})$ belongs to a certain parametric family, then the perturbation $\tilde{\boldsymbol{\omega}}$ in (25) is an appropriate perturbation. In particular, the maximum value of $\mathrm{FI}_{f,\mathbf{h}}$ is the score test statistic for testing the hypothesis $H_0 : \boldsymbol{\omega} = \mathbf{0}$.*

**4. Application to linear mixed models.** We consider data that are composed of a response $y_{ij}$ and a $q_1 \times 1$ covariate vector $\mathbf{x}_{ij}$ for $j = 1, \ldots, m_i$ within clusters $i = 1, \ldots, n$. We define the linear mixed models as

$$\mathbf{y}_i = \mathbf{x}_i \boldsymbol{\beta} + \boldsymbol{\epsilon}_i, \tag{26}$$

where $\mathbf{y}_i = (y_{i1}, \ldots, y_{im_i})^T$ is an $m_i \times 1$ vector, $\mathbf{x}_i^T = [\mathbf{x}_{i1}, \ldots, \mathbf{x}_{im_i}]$, $\boldsymbol{\beta}$ is a $q_1 \times 1$ vector of unknown parameters and $\boldsymbol{\epsilon}_i$ is normally distributed with mean zero and covariance matrix $\Sigma_i = \Sigma_i(\boldsymbol{\xi})$, in which $\boldsymbol{\xi}$ is a $q_2 \times 1$ vector. Thus, $\boldsymbol{\theta}^T = (\boldsymbol{\beta}^T, \boldsymbol{\xi}^T)$ is a $q \times 1$ vector, where $q = q_1 + q_2$.

For the linear mixed model, because the estimates of $\boldsymbol{\theta}$ (e.g., maximum likelihood estimates) may heavily depend on a small portion of the data or even one observation (or cluster), it is important to detect both influential clusters and influential individual observations. However, at either the subject or individual level, we cannot distinguish between influence due to the specific cluster characteristics and influence due to the characteristics of specific observations within a cluster. For further discussion on these issues, see Ouwens, Tan and Berger [21], Zhu and Lee [37, 38], Zhu, He and Fung [41] and Fung, Zhu, Wei and He [12], among many others.

The likelihood displacement function (Cook [7]) will be used throughout this section. We define $L(\boldsymbol{\theta}) = \log p(\mathbf{Y}|\boldsymbol{\theta})$ and $L(\boldsymbol{\theta}|\boldsymbol{\omega}) = \log p(\mathbf{Y}|\boldsymbol{\theta}, \boldsymbol{\omega})$. Let $\hat{\boldsymbol{\theta}}$ and $\hat{\boldsymbol{\theta}}_{\omega}$ be the maximum likelihood estimates of $L(\boldsymbol{\theta})$ and $L(\boldsymbol{\theta}|\boldsymbol{\omega})$, respectively; the likelihood displacement function is given by $LD(\boldsymbol{\omega}) = 2[L(\hat{\boldsymbol{\theta}}) - L(\hat{\boldsymbol{\theta}}_{\omega})]$. It can be shown that $H_{LD} = 2\Delta^T(-\ddot{L})^{-1}\Delta$, where $\Delta$ is a $q \times p$ matrix with elements $\partial^2 L(\boldsymbol{\theta}|\boldsymbol{\omega})/\partial\theta_i \, \partial\omega_j$ and $\ddot{L}$ is a $q \times q$ Hessian matrix with elements $\ddot{L}_{ij} = \partial^2 L(\boldsymbol{\theta})/\partial\theta_i \, \partial\theta_j$ evaluated at $\hat{\boldsymbol{\theta}}$ and $\boldsymbol{\omega}^0$. We calculate the Hessian matrix $-\ddot{L}$ as

$$-\ddot{L} \approx \sum_{i=1}^{n} \begin{bmatrix} \mathbf{x}_i^T \Sigma_i^{-1} \mathbf{x}_i & \mathbf{0} \\ \mathbf{0} & 0.5\partial_{\xi}\Sigma_i(\boldsymbol{\xi})(\Sigma_i^{-1} \otimes \Sigma_i^{-1})\partial_{\xi}\Sigma_i(\boldsymbol{\xi})^T \end{bmatrix},$$

where $\partial_{\xi}\Sigma_i(\boldsymbol{\xi}) = \partial \operatorname{vec}(\Sigma_i(\boldsymbol{\xi}))/\partial\xi$ is a $q_2 \times m_i^2$ matrix, in which we define

$$\operatorname{vec}(\mathbf{Z}) = (z_{11}, \ldots, z_{1m_i}, \ldots, z_{m_i 1}, \ldots, z_{m_i m_i})^T$$

for any $m_i \times m_i$ matrix $\mathbf{Z} = (z_{ij})$. We calculate the geometrical quantities of the perturbation manifold, the $\Delta = \partial^2 L/\partial\theta \, \partial\boldsymbol{\omega}$ matrix and the influence measures below.

4.1. *Perturbation of individual covariance matrix.* We consider the perturbation of the individual covariance matrix by assuming that

$$(27) \qquad \mathrm{Cov}(\mathbf{y}_i) = \omega_i^{-1} \Sigma_i \qquad \text{for all } i = 1, \dots, n.$$

Thus, $\boldsymbol{\omega}^0 = \mathbf{1}_n$ and $p = n$. For the perturbed model, both $L(\boldsymbol{\theta}|\boldsymbol{\omega})$ and $\ell(\boldsymbol{\omega}|\mathbf{Y}, \boldsymbol{\theta})$ equal

$$-\tfrac{1}{2} \sum_{i=1}^{n} \log |\Sigma_i| + \tfrac{1}{2} \sum_{i=1}^{n} n_i \log \omega_i - \tfrac{1}{2} \sum_{i=1}^{n} \omega_i (\mathbf{y}_i - \mathbf{x}_i \boldsymbol{\beta})^T \Sigma_i^{-1} (\mathbf{y}_i - \mathbf{x}_i \boldsymbol{\beta}).$$

After some calculations, we have $g_{ij}(\boldsymbol{\omega}) = 0.5 m_i \omega_i^{-2} \delta_{ij}$, $T_{ijk}(\boldsymbol{\omega}) = -m_i \omega_i^{-3} \times \delta_{ij} \delta_{ik}$ and $\Gamma_{ijk}^\alpha(\boldsymbol{\omega}) = 0.5(1-\alpha) T_{ijk}(\boldsymbol{\omega})$. At $\boldsymbol{\omega}^0$, $G = \mathrm{diag}(0.5 m_1, \dots, 0.5 m_n)$ indicates that the amount of perturbation introduced by $\omega_i$ is proportional to $m_i$, the number of observations in the $i$th cluster. Thus, $\boldsymbol{\omega}$ is an appropriate perturbation if and only if $m_1 = \dots = m_n$. Although $m_i$ may not be homogeneous, we can always consider an appropriate perturbation $\tilde{\boldsymbol{\omega}}$ in (5) such that $G(\tilde{\boldsymbol{\omega}})|_{\boldsymbol{\omega}^0} = \mathbf{I}_n$ and

$$(28) \qquad \mathrm{Cov}(\mathbf{y}_i) = [1 + (\tilde{\omega}_i - 1)/\sqrt{0.5 m_i}] \Sigma_i \qquad \text{for } i = 1, \dots, n.$$

For the appropriate perturbation (28), we get

$$\Delta_i = \frac{\partial^2 L(\boldsymbol{\theta}|\tilde{\boldsymbol{\omega}})}{\partial \boldsymbol{\theta} \, \partial \tilde{\omega}_i} = \sqrt{0.5} m_i^{-1/2} (2 \mathbf{e}_i^T \Sigma_i^{-1} \mathbf{x}_i, [\partial_\xi \Sigma_i (\Sigma_i^{-1} \otimes \Sigma_i^{-1}) \mathrm{vec}(\mathbf{e}_i \mathbf{e}_i^T)]^T)^T,$$

where $\mathbf{e}_i = \mathbf{y}_i - \mathbf{x}_i \boldsymbol{\beta}$. For simplicity, let $\boldsymbol{\beta}$ be the parameter of interest. It can be shown (Cook [7]) that

$$C_{\mathbf{E}_i} = \mathrm{SI}_{LD, \mathbf{E}_i} = 2 m_i^{-1} \mathbf{e}_i^T \Sigma_i^{-1/2} P_{ii} \Sigma_i^{-1/2} \mathbf{e}_i,$$

where $P_{ii} = \Sigma_i^{-1/2} \mathbf{x}_i (\sum_{j=1}^n \mathbf{x}_j^T \Sigma_j^{-1} \mathbf{x}_j)^{-1} \mathbf{x}_i^T \Sigma_i^{-1/2}$.

4.2. *Perturbation of responses.*

4.2.1. *Scheme one.* We consider the perturbation

$$(29) \qquad \mathbf{y}_i(\boldsymbol{\omega}) = \mathbf{y}_i + \omega_i \mathbf{1}_{m_i}.$$

Thus, $\boldsymbol{\omega}^0 = \mathbf{0}_n$ and $p = n$, where $\mathbf{0}_n$ represents an $n \times 1$ vector with all elements equal to 0. For the perturbed model, both $L(\boldsymbol{\theta}|\boldsymbol{\omega})$ and $\ell(\boldsymbol{\omega}|\mathbf{Y}, \boldsymbol{\theta})$ equal

$$-\tfrac{1}{2} \sum_{i=1}^{n} \log |\Sigma_i| - \tfrac{1}{2} \sum_{i=1}^{n} (\mathbf{y}_i + \omega_i \mathbf{1}_{m_i} - \mathbf{x}_i \boldsymbol{\beta})^T \Sigma_i^{-1} (\mathbf{y}_i + \omega_i \mathbf{1}_{m_i} - \mathbf{x}_i \boldsymbol{\beta}).$$

After some calculations, we have $g_{ij}(\boldsymbol{\omega}) = \mathbf{1}_{m_i}^T \Sigma_i^{-1} \mathbf{1}_{m_i} \delta_{ij}$, $T_{ijk}(\boldsymbol{\omega}) = 0$ and $\Gamma_{ijk}^\alpha(\boldsymbol{\omega}) = 0$. In this case, $G = \mathrm{diag}(\mathbf{1}_{m_1}^T \Sigma_1^{-1} \mathbf{1}_{m_1}, \dots, \mathbf{1}_{m_n}^T \Sigma_n^{-1} \mathbf{1}_{m_n})$ and the $i$th diagonal element of $G$ also depends on the number of observations in the $i$th cluster. This perturbation manifold is $\alpha$-flat for any $\alpha \in R^1$ and $\boldsymbol{\omega}(t) = t\mathbf{h}$ is a geodesic with respect to $\Gamma_{ijk}^\alpha$ for any $\alpha$. However, $\boldsymbol{\omega}$ is an appropriate perturbation if and only if $\mathbf{1}_{m_1}^T \Sigma_1^{-1} \mathbf{1}_{m_1} = \cdots = \mathbf{1}_{m_n}^T \Sigma_n^{-1} \mathbf{1}_{m_n}$. Therefore, $\boldsymbol{\omega}$ may be not an appropriate perturbation, but we can always consider an appropriate perturbation $\tilde{\boldsymbol{\omega}}$ in (5) such that $G(\tilde{\boldsymbol{\omega}})|_{\boldsymbol{\omega}^0} = \mathbf{I}_n$ and

$$(30) \qquad \mathbf{y}_i(\tilde{\boldsymbol{\omega}}) = \mathbf{y}_i + \omega_i \mathbf{1}_{m_i} / \sqrt{\mathbf{1}_{m_i}^T \Sigma_i^{-1} \mathbf{1}_{m_i}} \qquad \text{for } i = 1, \ldots, n.$$

For the perturbation $\tilde{\boldsymbol{\omega}}$ in (30), we have

$$\Delta_i = (\mathbf{1}_{m_i}^T \Sigma_i^{-1} \mathbf{x}_i, [\partial_\xi \Sigma_i (\Sigma_i^{-1} \otimes \Sigma_i^{-1}) \operatorname{vec}(\mathbf{1}_{m_i} \mathbf{e}_i^T)]^T)^T.$$

Let $\boldsymbol{\xi}$ be the parameter of interest. It can be shown that

$$\begin{aligned} C_{\mathbf{E}_i} = {}& 4(\mathbf{1}_{m_i}^T \Sigma_i^{-1} \mathbf{1}_{m_i})^{-1} \\ & \times \operatorname{vec}(\mathbf{1}_{m_i} \mathbf{e}_i^T)^T (\Sigma_i^{-1/2} \otimes \Sigma_i^{-1/2}) Q_{ii} (\Sigma_i^{-1/2} \otimes \Sigma_i^{-1/2}) \operatorname{vec}(\mathbf{1}_{m_i} \mathbf{e}_i^T), \end{aligned}$$

where $Q_{ii}$ is defined as

$$\begin{aligned} (\Sigma_i^{-1/2} & \otimes \Sigma_i^{-1/2}) (\partial_\xi \Sigma_i)^T \left[ \sum_{j=1}^n (\partial_\xi \Sigma_j)(\Sigma_j^{-1} \otimes \Sigma_j^{-1})(\partial_\xi \Sigma_j)^T \right]^{-1} \\ & \times (\partial_\xi \Sigma_i)(\Sigma_i^{-1/2} \otimes \Sigma_i^{-1/2}). \end{aligned}$$

### 4.2.2. *Scheme two.* We consider the mean shift perturbation model

$$(31) \qquad \mathbf{y}_i(\boldsymbol{\omega}) = \mathbf{y}_i + \boldsymbol{\omega}_i,$$

where $\boldsymbol{\omega}_i = (\omega_{i1}, \ldots, \omega_{im_i})^T$. Thus, $\boldsymbol{\omega}^{0T} = (\mathbf{0}_{m_1}^T, \ldots, \mathbf{0}_{m_n}^T)$ and $p = \sum_{i=1}^n m_i$. For the perturbed model, both $L(\boldsymbol{\theta}|\boldsymbol{\omega})$ and $\ell(\boldsymbol{\omega}|\mathbf{Y}, \boldsymbol{\theta})$ equal

$$-\frac{1}{2} \sum_{i=1}^n \log|\Sigma_i| - \frac{1}{2} \sum_{i=1}^n (\mathbf{y}_i + \boldsymbol{\omega}_i - \mathbf{x}_i \boldsymbol{\beta})^T \Sigma_i^{-1} (\mathbf{y}_i + \boldsymbol{\omega}_i - \mathbf{x}_i \boldsymbol{\beta}).$$

After some calculations, we have $g_{ij}(\boldsymbol{\omega}) = \Sigma_i^{-1} \delta_{ij}$, $T_{ijk}(\boldsymbol{\omega}) = 0$ and $\Gamma_{ijk}^\alpha(\boldsymbol{\omega}) = 0$, where $i$, $j$ and $k$ vary from 1 to $n$. The structure of the metric tensor $G(\boldsymbol{\omega}) = (g_{ij}(\boldsymbol{\omega}))$ indicates that the perturbations $\boldsymbol{\omega}_i$ in different clusters are orthogonal to each other, whereas the components $\{\boldsymbol{\omega}_{il} : l = 1, \ldots, m_i\}$ of $\boldsymbol{\omega}_i$ are associated with each other. This perturbation manifold is also $\alpha$-flat for any $\alpha \in R^1$ and $\boldsymbol{\omega}(t) = t\mathbf{h}$ is a geodesic with respect to $\Gamma_{ijk}^\alpha(\omega)$ for any $\alpha$. However, $\boldsymbol{\omega}$ is not an appropriate perturbation, because $G(\boldsymbol{\omega})|_{\boldsymbol{\omega}^0}$ does not have the form $c\mathbf{I}_M$, where $M = \sum_{i=1}^n m_i$. Therefore, we consider an appropriate perturbation $\tilde{\boldsymbol{\omega}}$ in (5) with $c = 1$ such that $G(\tilde{\boldsymbol{\omega}})|_{\boldsymbol{\omega}^0} = \mathbf{I}_n$ and

$$(32) \qquad \mathbf{y}_i(\tilde{\boldsymbol{\omega}}) = \mathbf{y}_i + \Sigma_i^{-1/2} \tilde{\boldsymbol{\omega}}_i \qquad \text{for } i = 1, \ldots, n.$$

For the perturbation $\tilde{\boldsymbol{\omega}}$ in (32), we have $\Delta_i = [\Sigma_i^{-1}\mathbf{x}_i, (\mathbf{e}_i^T\Sigma_i^{-1}\otimes\Sigma_i^{-1})\partial_\xi\Sigma_i^T]^T$. Let $\boldsymbol{\xi}$ be the parameter of interest. It can be shown that (Cook [7])

$$\mathbf{R}_i = 4\Sigma_i^{1/2}(\mathbf{e}_i^T\Sigma_i^{-1/2}\otimes\Sigma_i^{-1/2})Q_{ii}(\Sigma_i^{-1/2}\mathbf{e}_i\otimes\Sigma_i^{-1/2})\Sigma_i^{1/2}$$

is the submatrix of $H_{LD}$ corresponding to the $i$th perturbation vector $\tilde{\boldsymbol{\omega}}_i$. Therefore, $C_{\mathbf{E}_{i,l}} = \mathrm{SI}_{LD,\mathbf{E}_{i,l}}$, which corresponds to the $l$th diagonal element of $\mathbf{R}_i$, where $\mathbf{E}_{i,l}$ is a $p\times 1$ vector with a 1 at the $(\sum_{k=1}^{i-1}m_k+l)$th element and 0 otherwise.

### 4.3. Yale infant growth data.

The Yale infant growth data were collected to study whether cocaine exposure during pregnancy may lead to the maltreatment of infants after birth, such as physical and sexual abuse. A total of 298 children were recruited from two subject groups (cocaine exposure group and unexposed group). The key feature of this dataset is that different children had different numbers and patterns of visits during the study period. We refer to Wasserman and Leventhal [29] and Stier et al. [24] for a detailed description of the study design and data collection. Recently, Zhang [34, 35] developed multivariate adaptive splines for the analysis of longitudinal data (MASAL) to analyze the Yale infant growth data. The importance of our reanalysis here is to develop a local influence approach for the MASAL model.

For the Yale infant growth data, Zhang [35] selected the MASAL model

$$y_{ij} = \mathbf{x}_{ij}^T\boldsymbol{\beta} + \epsilon_{i,j},$$

where $\mathbf{x}_{ij} = (1, d, (d-120)^+, (d-200)^+, (g_a-28)^+, d(g_a-28)^+, (d-60)^+(g_a-28)^+, (d-490)^+(g_a-28)^+, sd, s(d-120)^+)^T$, in which $d$ and $g_a$ are the age of visit and gestational age, respectively, and $s$ is the indicator for gender, with one indicating a girl and zero indicating a boy. In addition, we assume that $\boldsymbol{\epsilon}_i = (\epsilon_{i1},\ldots,\epsilon_{im_i})^T \sim N[\mathbf{0}, \Sigma_i(\boldsymbol{\xi})]$ and $\Sigma_i(\boldsymbol{\xi})$ is determined by the variance and autocorrelation parameters, which are, respectively, given as

$$V(d) = \exp(\xi_0 + \xi_1 d + \xi_2 d^2 + \xi_3 d^3) \quad \text{and} \quad \rho(l) = \xi_4 + \xi_5 l,$$

where $l$ is the lag between two visits. For simplicity, we assume that all knots are given so that the MASAL model reduces to the linear mixed model (26). The total number of data points is $\sum_{i=1}^n m_i = 3176$ and the total number of clusters is $n = 298$. The estimated parameters are

$$\boldsymbol{\beta}^T = (0.744, 0.029, -0.0092, -0.0059, 0.204,$$
$$0.0005, -0.0007, -0.0009, -0.0026, 0.0022)$$

and

$$\boldsymbol{\xi} = (-0.53, 0.0064, -1.9\times 10^{-5}, 2.1\times 10^{-8}, 0.929, -0.0013)^T$$

(Zhang [35]).

We calculated the local influence measures for the three perturbations discussed in Sections 4.1 and 4.2 and present the main findings of the local influence approach in Figures 1–3.

For the perturbation of the individual covariance matrix, the quantity $SI_{\mathbf{E}_i}$ for perturbation (28) reveals four influential subjects $\{141, 246, 269, 285\}$ (Figure 1), whereas subjects $\{141, 246\}$ do not stand out as influential using normal curvature for perturbation (27). Because $m_{141} = 4$ and $m_{246} = 5$ are much smaller than the average number of observations $3176/298 = 10.6$, a relatively large normal curvature from subjects $\{141, 246\}$ represents a large effect. A closer inspection of the data (not presented here) shows that the raw and fitted curves for these four subjects differ substantially, especially at the last observation.

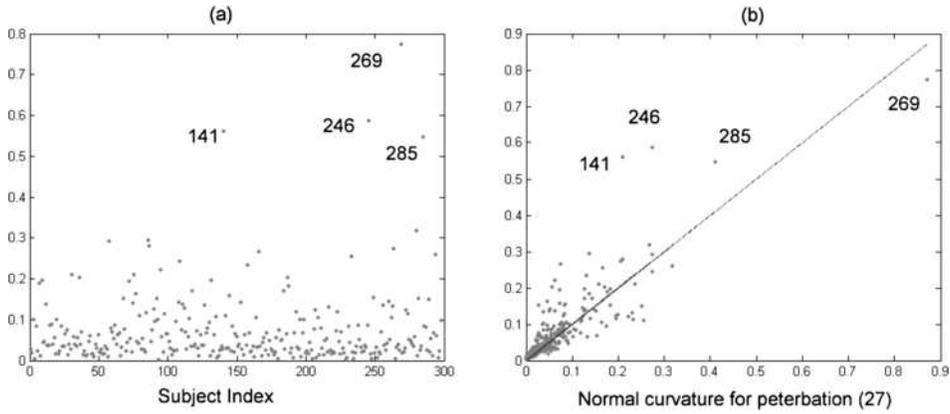

Fig. 1. *Yale infant growth data:* (a) *index plot of* $3176 \times SI_{\mathbf{E}_i}/298$ *for the appropriate perturbation* (28); (b) $C_{\mathbf{E}_i}$ *for perturbation* (27) *against* $3176 \times SI_{\mathbf{E}_i}/298$ *for the appropriate perturbation* (28).

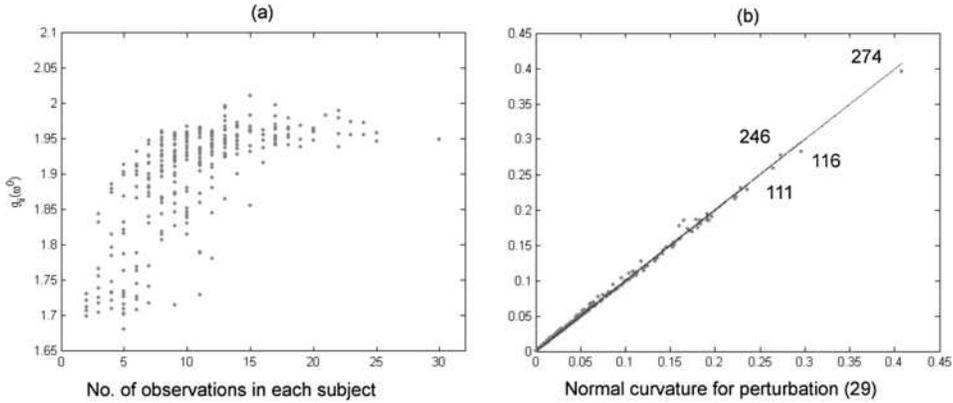

Fig. 2. *Yale infant growth dataset:* (a) *metric tensor* $g_{ii}(\boldsymbol{\omega}^0)$ *for perturbation* (29) *and the number of observations for each subject;* (b) $C_{\mathbf{E}_i}$ *for perturbation* (29) *and* $1.9 \times SI_{\mathbf{E}_i}$ *for the appropriate perturbation* (30).

The metric tensor $g_{ii}(\omega^0)$ for perturbation (29) is positively correlated with the number of observations in each subject [Figure 2(a)]. Furthermore, because the variability of $g_{ii}(\omega^0)$ is relatively small, normal curvatures under perturbation (29) are close to second-order influence measures [Figure 2(b)]. Both $C_{\mathbf{E}_i}$ and $\mathrm{SI}_{\mathbf{E}_i}$ reveal four influential subjects $\{111, 116, 246, 274\}$ [Figure 2(b)]. Furthermore, for perturbation (32), $\mathrm{SI}_{\mathbf{E}_i}$ suggest that $(3, 7)$, $(24, 7)$, $(24, 8)$, $(227, 11)$, $(290, 12)$, $(290, 13)$ and $(109, 12)$ are seven influential observations, where for each $(i, l)$, $i$ denotes the subject number and $l$ denotes the observation number [Figure 3(a)–(b)].

## 5. Conclusion.

We have introduced a local influence method to assess minor perturbations to a statistical model. Our method extends the previous local influence method (Cook [7]) in several aspects. First, we propose to use the metric tensor of a perturbation manifold to select an appropriate perturbation to a model. The major advantage of using an appropriate perturbation is that it leads to a nice interpretation of the effect of all elements of a perturbation vector on a statistical model. We have shown in Sections 3.1–3.3 that most of the perturbation schemes considered in Cook's [7] examples are appropriate. However, we have also shown in several examples, such as linear mixed models and testing parametric families, that some commonly used perturbations may not yield an appropriate perturbation; see Sections 3.4 and 4. Second, we have developed influence measures with nice geometrical interpretations for smooth objective functions at any point. The influence measures proposed here avoid the previous drawback that the normal curvature is not well defined for some objective functions at points with a nonzero first derivative. In addition, the proposed second-order influence measures reduce to normal curvatures for the likelihood displacement function (Cook [7]). Third, we have established a connection between the score test statistic and the FI measures; see Section 3.4. Finally, we have also examined a number of models to highlight the importance of choosing an appropriate perturbation and the broad spectrum of applications of this local influence method.

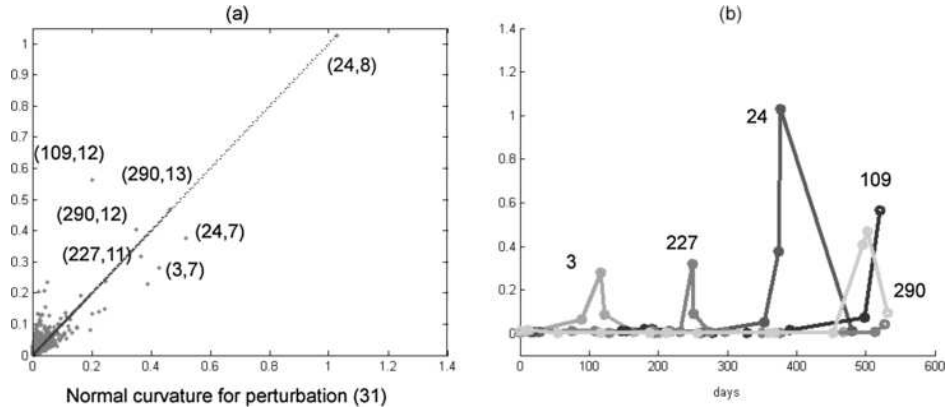

FIG. 3. *Yale infant growth data:* (a) *rescaled* $C_{\mathbf{E}_{i,l}}$ *for perturbation* (31) *and* $\mathrm{SI}_{\mathbf{E}_{i,l}}$ *for the appropriate perturbation* (32)*;* (b) *connected curves of* $\mathrm{SI}_{\mathbf{E}_{i,l}}$ *of five influential subjects* $\{3, 24, 109, 227, 290\}$ *for the appropriate perturbation* (32)*.*

Many issues still merit further research. One major issue is calculation of the influence measures and metric tensor under different situations, such as measurement error models (Carroll, Ruppert and Stefanski [4]; Zhong, Wei and Fung [36]), generalized linear models with missing data (Ibrahim, Chen, Lipsitz and Herring [13]), partially linear models (Zhu, He and Fung [41]) and structural equation models (Yuan and Bentler [33]). Another major issue is to establish relationships between the influence measures and other influence diagnostics, such as case-deletion measures and leverage (Cook and Weisberg [8] and Wei, Hu and Fung [30]). It is also important to develop appropriate influence diagnostics for detecting influential clusters in longitudinal data by taking into account the number of observations in each cluster and models used to fit the longitudinal data. However, the influence diagnostics calculated in PROC MIXED of SAS 9.1 (SAS Institute Inc., Cary, NC) do not take into account the number of observations in each cluster; therefore, these influence diagnostics may give misleading results. We expect that the metric tensor of the perturbation manifold will play a critical role in this new development.

**Acknowledgments.** We thank the Editor Morris Eaton, the Associate Editor and two anonymous referees for valuable suggestions, which greatly helped to improve our presentation. Reprints can be requested via e-mail: heping.zhang@yale.edu or htzhu@email.unc.edu.

H. ZHU
DEPARTMENT OF BIOSTATISTICS
    AND BIOMEDICAL RESEARCH IMAGING CENTER
UNIVERSITY OF NORTH CAROLINA AT CHAPEL HILL
CHAPEL HILL, NORTH CAROLINA 27599-7420
USA
E-MAIL: htzhu@email.unc.edu

S. LEE
DEPARTMENT OF STATISTICS
THE CHINESE UNIVERSITY OF HONG KONG
HONG KONG
CHINA
E-MAIL: sylee@sparc2.sta.cuhk.edu.hk

J. G. IBRAHIM
DEPARTMENT OF BIOSTATISTICS
UNIVERSITY OF NORTH CAROLINA
    AT CHAPEL HILL
CHAPEL HILL, NORTH CAROLINA 27599-7420
USA
E-MAIL: ibrahim@bios.unc.edu

H. ZHANG
DEPARTMENT OF EPIDEMIOLOGY
    AND PUBLIC HEALTH
YALE UNIVERSITY SCHOOL OF MEDICINE
NEW HAVEN, CONNECTICUT 06520-8034
USA
AND
JIANGXI NORMAL UNIVERSITY
NANCHANG
CHINA
E-MAIL: heping.zhang@yale.edu